\documentclass[12pt,leqno]{article}
\usepackage{latexsym,amsfonts,epsfig}
\usepackage{amsmath,amssymb,amsthm}
\usepackage[notref, notcite]{}
\usepackage[colorlinks,linkcolor=red,citecolor=blue,urlcolor=blue]{hyperref}
\usepackage[margin=0.5in]{geometry}
\input epsf
\newtheorem{theorem}{Theorem}[section]
\newtheorem{lemma}{Lemma}[section]
\newtheorem{corollary}{Corollary}[section]
\newtheorem{remark}{Remark}[section]

\begin{document}
\title{Sharp Inequalities of Bienaym\'e-Chebyshev and Gau\ss \, Type for Possibly Asymmetric Intervals around the Mean}
\author{  Roxana A. Ion (ASML Veldhoven) \\
 Chris A.J. Klaassen (University of Amsterdam) \\
           Edwin R. van den Heuvel (Technical University Eindhoven)}

\maketitle



\begin{abstract}
Gau\ss \, (1823) proved a sharp upper bound on the probability that a random variable falls outside a symmetric interval around zero when its distribution is unimodal with mode at zero.
For the class of all distributions with mean at zero, Bienaym\'e (1853) and Chebyshev (1867) independently provided another, simpler sharp upper bound on this probability.
For the same class of distributions, Cantelli (1928) obtained a strict upper bound for intervals that are a half line.
We extend these results to arbitrary intervals for six classes of distributions, namely the general class of `distributions', the class of `symmetric distributions', of `concave distributions', of `unimodal distributions', of `unimodal distributions with coinciding mode and mean', and of `symmetric unimodal distributions'.
For some of the known inequalities, such as the Gau\ss \, inequality, an alternative proof is given.
\end{abstract}

\section{Introduction}\label{introduction}

Let $X$  be a random variable with finite mean $\mu$ and finite, positive variance $\sigma^2$.
We are interested in sharp upper bounds on the tails of the distribution of $X$.
Without loss of generality and in order to simplify notation we will restrict attention to the standardized version of $X$, which is denoted by $Z=(X-\mu)/\sigma.$
Hence, we are interested in sharp upper bounds on the probability that a standardized random variable falls outside an arbitrary interval containing the value zero, i.e. sharp upper bounds on
\begin{equation}\label{rint1}
P(Z\leq -u\ \  \mbox{\rm or}\ \  Z \geq v), \quad u>0,\,\, v>0.
\end{equation}
We may assume $0<v \leq u$.
We also study the one-sided probability $P(Z \geq v)$, which corresponds to (\ref{rint1}) with $u=\infty$.

Of course, upper bounds on probability $(\ref{rint1})$ depend on the class of distributions assumed for $Z$.
We shall focus on broad, nonparametric classes of distributions.
The corresponding upper bounds are often applied in theoretical considerations and proofs, typically with $u=v$.
However, as \cite{Rougier}, p.861, argue, they have practical value as well when one does not want strong model assumptions.
Let us mention two examples, one with $u$ equal to infinity and the other one with general $u$ and $v$.

The authors of \cite{Clarkson} study the Receiver Operating Characteristic (ROC) curve, which they define as the power of a simple versus simple Neyman-Pearson hypothesis test viewed as a function of the significance level, more precisely
\begin{eqnarray*}
\alpha \mapsto 1 - G(F^{-1}(1-\alpha)),\quad \alpha \in [0,1],
\end{eqnarray*}
where $F$ and $G$ are the distribution functions of the test statistic under the null hypothesis and the alternative, respectively.
They are interested in the Area Under Curve (AUC), which equals
\begin{eqnarray*}\label{AUC}
\int_0^1 \left(1 - G(F^{-1}(1-\alpha))\right) d\alpha = P(X \leq Y)
\end{eqnarray*}
with $X$ and $Y$ independent random variables with distribution functions $F$ and $G$ respectively.
For Gaussian $X$ and $Y$ with variance 1 and $E(Y-X)= 2\sqrt{6}$ the AUC equals $\Phi(2\sqrt{3})\approx 0.99973$ (and not 0.99966 as at page 468 of \cite{Clarkson}).
If Gaussianity is weakened to unimodality of $X$ and $Y$ with at least one of them being strongly unimodal, then \cite{Clarkson} obtain 1/2 as a lower bound to the AUC; see their page 468, Corollary 6 and Remark 7.
However, our bound (\ref{unimodalCantellia}) from Theorem \ref{unimodalCantelli} with $v=2\sqrt{3}$ yields $113/117 \approx 0.96581$ as a lower bound.
This shows that at least some of the bounds presented here are strong enough to be valuable in applications.

Sharp upper bounds on probability (\ref{rint1}) are also relevant in statistical process control.
Indeed, attainable upper bounds on (\ref{rint1}) determine the maximum risk for producing products outside specification limits, as the distribution of the quality characteristic $X$ is typically not (perfectly) centered within the specification limits.
This kind of practical concerns was our main motivation for studying upper bounds on (\ref{rint1}); see for instance
\cite{Heuvel} and the PhD thesis \cite{Ion}.

Traditionally, focus has been on upper bounds on the probability in (\ref{rint1}) for symmetric intervals with $u=v.$
The Bienaym\'e-Chebyshev inequality, which may be stated and proved by
\begin{equation}\label{BCi}
 P(\lvert Z \rvert  \geq v) = E \left( {\bf 1}_{[Z^2/v^2 \geq 1]} \right) \leq  E \left( \frac{Z^2}{v^2}{\bf 1}_{[Z^2/v^2 \geq 1]} \right)
 \leq  E \left( \frac{Z^2}{v^2} \right) = \frac{1}{v^2},
\end{equation}
was first obtained by Bienaym\'e \cite{Bienayme}, cf. page 171 of the 1867 reprint.
Hidden in a discourse about convergence of the least squares method, he formulates his proof as follows.
On the one hand $E(Z^2{\bf 1}_{[Z^2/v^2 \geq 1]})$ is a fraction $f$ of $EZ^2$.
On the other hand, as $v^2 {\bf 1}_{[Z^2/v^2 \geq 1]}$ is the smallest value $Z^2 {\bf 1}_{[Z^2/v^2 \geq 1]}$ can take, there exists a (positive) constant $\theta \leq 1$ such that $E(Z^2{\bf 1}_{[Z^2/v^2 \geq 1]})$ equals $P(\lvert Z \rvert   \geq v) v^2/\theta$.
Bienaym\'e concludes $P(\lvert Z \rvert   \geq v) = \theta f / v^2$ with $\theta f \leq 1$.
In his proof of the law of large numbers, actually Chebyshev \cite{Chebyshev} followed the line of (\ref{BCi}), but he argued more straightforwardly by noting that $Z^2/v^2 \geq {\bf 1}_{[Z^2/v^2 \geq 1]}$ holds.
See also \cite{Heyde}, page 682.
Note that in (\ref{BCi}) equality holds if and only if $\lvert Z \rvert  $ takes on the values 0 and $v$ only.

\begin{table}\label{table1}
{\centering
\begin{tabular}{ |c|c|c|c|c|}
\hline
Class & & & \\
of & $P(Z\geq v)$ & $P(\lvert Z \rvert  \geq v)$ &
$P(Z\leq -u,Z\geq v)$ \\
distributions & & & \\
\hline
& $\frac{1}{1+v^2}$ & $\frac{1}{v^{2}}$ & $\frac{4+(u-v)^2}{(u+v)^2}$ \\
All & Cantelli (1928)  & Bienaym\'e (1853) & \\
& Theorem \ref{Cantelli} & Chebyshev (1867) &  Theorem \ref{general}\\
\hline
& $\frac{1}{2v^2}$ & $\frac{1}{v^2}$ & $\max\{\frac{1}{u^2},\frac{1}{2v^2}\}$  \\
Symmetric & Chebyshev (1867) & Chebyshev (1867) &  $1 \leq v \leq u$ \\
& Theorem \ref{symmetricChebyshev} & & Theorem \ref{symmetricgeneral} \\
\hline
Concave & $\frac{4}{9}\frac{1}{v^{2}}$ && \\
on $[0,\infty)$  & & irrelevant & irrelevant \\
 $EZ>0$ &  Theorem \ref{onesidedGauss} && \\
\hline
& $\frac{4}{9}\frac{1}{1+v^{2}}$ & $\frac 49 \frac{1}{v^{2}}$ & $\frac 49 \frac{4+(u-v)^2}{(u+v)^2}$ \\
Unimodal & & Vyso$\check{\rm c}$anski$\breve{\rm i}$- &  [Vysochanski$\breve{\rm i}$- \\
         & & Petunin (1980) & Petunin (1983)] \\
& Theorem \ref{unimodalCantelli} & Corollary \ref{absoluteGauss} & Theorem \ref{generalGauss} \\
\hline
Unimodal & (\ref{boundonesidedGauss}) & $\frac{4}{9}\frac{1}{v^{2}}$ & (\ref{boundtwosidedGauss})--(\ref{valuec}) \\
\& &&  Gau\ss\,(1823) &  \\
mode = mean & Theorem \ref{onesidedGauss2} & Theorem \ref{Gauss} & Theorem \ref{twosidedGauss2} \\
\hline
Symmetric & $\frac{4}{9}\frac{1}{2v^{2}}$ & $\frac{4}{9}\frac{1}{v^{2}}$ & $\frac 49 \frac 4{(u+v)^2}$ \\
\& & Camp (1922) & Gau\ss \, (1823) & Semenikhin (2019) \\
unimodal & Theorem \ref{Gausslike} & Theorem \ref{Gauss}  & Theorem \ref{symmetricunimodaltheorem} \\
\hline
\end{tabular}\par}
\
\caption{Overview of the inequalities of Bienaym\'e-Chebyshev and Gau\ss \, type for standardized random variables $Z$ (except for the concave case).}
\end{table}

Amazingly, Gau\ss \, \cite{Gauss}, translation \cite{Stewart}, had already proved the sharp inequality
\begin{equation}\label{Gauss0}
P(\lvert Y \rvert   \geq v) \leq \frac 49 \ \frac 1{v^2}
\end{equation}
for random variables $Y$ with a unimodal distribution with mode at zero and with $E(Y^2) = 1$.
The authors of \cite{Sellke2} describe the history of the Gau\ss \, inequality and refer to \cite{Pukelsheim} for three proofs, namely Gau\ss's one, a variation on it, and the one from exercise 4 on page 256 of \cite{Cramer}.
A detailed description of Gau\ss's proof is given also by \cite{Hooghiemstra}.
We present a fourth proof of the Gau\ss \, inequality in Section \ref{concave}.
However, as the mean is typically known or estimated in practice (and can be estimated more accurately than the mode), we shall focus on standardized random variables $Z$ with $EZ=0$.
For the asymmetric case with $u \neq v$ in (\ref{rint1}) this mean zero condition complicates results and proofs considerably; see Theorem \ref{twosidedGauss2}.

The upper bound in $(\ref{BCi})$ follows from the classical Markov inequality \\
$P(\lvert Z \rvert  \geq v) = E({\bf 1}_{[\lvert Z \rvert  ^r/v^r \geq 1]}) \leq E\lvert Z \rvert  ^r/v^r$ with $r=2.$
Generalizations of the Gau\ss \ inequality in Markov style have been presented by \cite{Camp}, \cite{Meidell}, and \cite{Theil}.
\cite{DasGupta} has exploited these Markov-Gau\ss \,-Camp-Meidell inequalities to derive properties of e.g. the zeta function.
\cite{Sellke1} and \cite{Sellke2} extend this line of research by determining sharp upper bounds on $P(\lvert Z \rvert  \geq v)$ in terms of
$Eg( \lvert Z \rvert )$ for unimodal random variables $Z$ and nondecreasing functions $g$ on $[0,\infty).$
\cite{Bickel} derives improvements of (\ref{BCi}) with $Z$ an average and of (\ref{Gauss0}) for symmetric unimodal densities with a bound on the derivative.

For symmetric distribution functions of $X$, an upper bound on the one-sided probability $P(Z\geq v)$  is  $1/(2v^2)$, which can be
obtained by the inequality in $(\ref{BCi})$, see Theorem \ref{symmetricChebyshev}.
Cantelli \cite{Cantelli} proved the upper bound $1/(1+v^2)$  on this one-sided probability for the class of all distribution functions of $X.$
\cite{Camp} provided an upper bound when the distribution function of $X$ is symmetric and unimodal.

In each of the subsequent Sections we shall discuss sharp inequalities for a specific class of distributions of the standardized random variable $Z$.
The only exception is in Section \ref{concave}, where we do not consider standardized random variables.
There we shall present a one-sided, sharp version of Gau\ss's inequality for concave distributions, more specifically, for distributions on the nonnegative half line with a mode at zero and a concave distribution function on this half line.
As an immediate consequence of this one-sided version, the original Gau\ss \, inequality is also presented in this Section.
Our proof is based on Khintchine's representation of unimodal densities, \cite{Khintchine}, and Jensen's inequality, \cite{Jensen}, as are our proofs of most results for unimodal distributions.
A formulation and proof of Khintchine's representation may be found in Lemma \ref{Khintchinelemma} in the Appendix.

The results on upper bounds for the probability in $(\ref{rint1})$ as discussed in this article, are summarized in Table \ref{table1}, where we mention only the most relevant cases of the inequalities.
The complete versions of the inequalities with all special cases and proofs may be found in the text.
Reference is made to the corresponding theorems that are proved in the present article.
Theorems \ref{symmetricChebyshev}, \ref{symmetricgeneral}, \ref{onesidedGauss}, \ref{unimodalCantelli}, \ref{generalGaussa}, \ref{generalGauss}, \ref{onesidedGauss2} and \ref{twosidedGauss2}, and Lemma \ref{upperboundconvex} seem to be new.
For known results like in Theorems \ref{general}, \ref{originalGauss}, \ref{Gausslike} and \ref{symmetricunimodaltheorem} alternative proofs are given.
It is shown that all inequalities are sharp as well, by constructing random variables satisfying the bounds.
It should be noted that the results on probability $(\ref{rint1})$ imply the results for upper bounds on the probability in $(\ref{BCi})$.

\section{All distributions}\label{all}

Here we discuss the inequalities from the first row of Table \ref{table1}.
The standardized version $Z=(X-\mu) / \sigma $ of the random variable $X$ has mean $0$ and variance $1$.
Apart from this standardization the distribution of $Z$ is arbitrary within this section.
We will start with the one-sided analogue of the Bienaym\'e-Chebyshev inequality, which is referred to as Cantelli's inequality.
It is given by formula (19) in \cite{Cantelli} and it reads as follows.

\begin{theorem}{\bf Cantelli's inequality}\label{Cantelli}\\
Let $Z$ be a random variable with mean 0 and variance 1.
For any $v\geq 0$ the inequality
\begin{equation}\label{cantelli}
P(Z\geq v)\leq \frac{1}{ 1+v^2}
\end{equation}
holds. This inequality is sharp and for positive $v$ equality is uniquely attained by
\begin{equation}
\label{c43} Z=\left\{\begin{array}{rcl}
v & & \frac 1{1+v^2}\\
& \mbox{with probability} & \\
-\frac 1v &  & \frac{v^2}{1+v^2}.\\
\end{array}\right.
\end{equation}
\end{theorem}
{\bf Proof}\\
According to the hint from Problem 1.5.5 of \cite{Billingsley}, the Bienaym\'e-Chebyshev inequality yields
\begin{equation}\label{c51}
P(Z\geq v) = P(1+vZ \geq 1+v^2) \leq E\left(\frac{1+vZ}{1+v^2}\right)^2 = \frac{1}{1+v^2}\,.
\end{equation}
In order to obtain equality in $(\ref{c51})$ the equality $$ {\bf 1}_{[Z\geq v]}=((1+vZ)/(1+v^2))^2 $$ has to hold almost surely, or equivalently, $Z$ has to have support $\{-1/v, v \}$, and hence has to satisfy $(\ref{c43})$.
\hfill$\Box$\\

The asymmetric two-sided analogue of the Bienaym\'e-Chebyshev inequality from (\ref{BCi}) will be presented in full detail.
Its simple proof is based on the Bienaym\'e-Chebyshev inequality itself.

\begin{theorem}\label{general}
Let $Z$ have mean 0 and variance 1, and assume $0<v \leq u$.
Then
\begin{equation}\label{fourplus}
P(Z\leq -u\ \  \mbox{\rm or}\ \  Z \geq v)\leq
\frac{4+(u-v)^2}{(u+v)^2}
\end{equation}
holds.
Under the additional condition $ u \leq 1/v$ the trivial inequality
\begin{equation}\label{r1412}
P(Z\leq -u\ \  \mbox{\rm or}\ \  Z\geq v)\leq 1
\end{equation}
holds with equality for
\begin{eqnarray}\label{r151}
Z=\left\{\begin{array}{rcl}
-\sqrt{\frac{u}{v}} & & \frac{v}{u+v} \\
 & \mbox{with probability} & \\
\sqrt{\frac{v}{u}} &  & \frac{u}{u+v}\,.\\
\end{array}\right.
\end{eqnarray}
In the case of
\begin{equation}\label{r1211}
\frac 1v \leq u \leq v+\frac 2v
\end{equation}
inequality (\ref{fourplus}) is sharp and equality holds for
\begin{eqnarray}\label{r1311}
\qquad Z=\left\{\begin{array}{rcl}
-u & & \frac{2+v(v-u)}{(u+v)^{2}} \\
\\
\frac {v-u}{2} & \mbox{with probability} & 1-\frac{4+(u-v)^2}{(u+v)^{2}} \\
\\
v &  & \frac{2+u(u-v)}{(u+v)^{2}}.\\
\end{array}\right.
\end{eqnarray}
In the case of $v+2/v \leq u$ the inequality
\begin{equation}\label{r181}
P(Z\leq -u\ \ \mbox{\rm or}\ \ Z\geq v)\leq \frac{1}{1+v^2}
\end{equation}
holds with equality for
\begin{eqnarray}\label{r191}
Z=\left\{\begin{array}{rcl}
            - \frac 1v & & \frac{v^2}{1+v^2} \\
                & \mbox{with probability} & \\
            v &  & \frac 1{1+v^2}\,.\\
\end{array}\right.
\end{eqnarray}
\end{theorem}
{\bf Proof}\\
By the Bienaym\'e-Chebyshev inequality we obtain
\begin{eqnarray*}
\lefteqn{P(Z \leq -u \ \ \mbox{\rm or} \ \  Z \geq v) = P \left( \mid Z - \tfrac 12 (v-u) \mid \geq \tfrac 12 (u+v) \right) } \\
&& = P\left( \mid 2Z - (v-u) \mid \geq u+v \right) \leq E\left( \frac{2Z-(v-u)}{u+v} \right)^2 = \frac{4+(v-u)^2}{(u+v)^2}.
\end{eqnarray*}
Straightforward computations show that $(\ref{r1311})$ yields a well-defined random variable provided that $(\ref{r1211})$ holds,
and that this random variable has mean zero and unit variance, and satisfies equality in $(\ref{fourplus})$.

Similarly, the random variable from $(\ref{r151})$ is well-defined and attains equality in the trivial inequality $(\ref{r1412})$ in
view of $-\sqrt {u/v} \leq -u$ iff $uv\leq 1$ iff $\sqrt {v/u} \geq v$.

For $u \geq v+2/v$ the probability in (\ref{fourplus}) is bounded by $P(Z\leq -v-2/v\ \ \mbox{or}\ \  Z\geq v)$
Applying inequality (\ref{fourplus}) with $u$ replaced by $v+2/v$ results in $(\ref{r181})$.
Finally, straightforward computations show that $(\ref{r191})$ defines a proper random variable with zero mean and unit variance that
attains equality in $(\ref{r181}).$
\hfill$\Box$\\

For $v$ fixed, the minimum of the right hand side of (\ref{fourplus}) over $u$ is attained at $u=v+2/v$ and equals $1/(1+v^2).$
Consequently, Cantelli's inequality (\ref{cantelli}) is a special case of (\ref{fourplus}); cf. (\ref{r181}) and (\ref{r191}).

Selberg \cite{Selberg} seems to be the first to have formulated a version of Theorem \ref{general}.
According to Ferentinos \cite{Ferentinos} his proof of (\ref{fourplus}) is less complicated than Selberg's, but it is still more complicated than ours due to the cumbersome notation.

Note that (\ref{fourplus}) with $u=v$ reduces to the famous Bienaym\'e-Chebyshev inequality.
The inequalities in this section are all based on the inequality of Bienaym\'e-Chebyshev itself, as is the first one of the next section.
However, this inequality does not always seem to be helpful if the class of distributions of $Z$ (or $X$) is restricted.

\section{Symmetric distributions}\label{symmetric}
%\label{symmetricdistributions}

The inequalities for symmetric distributions from the second row of Table \ref{table1} will be discussed in this section.
The symmetry implies $P(Z \geq v) = P(\lvert Z \rvert  \geq v)/2$ and hence the Bienaym\'e-Chebyshev inequality yields the following result.

\begin{theorem}\label{symmetricChebyshev}
Let $Z$ be symmetric with mean 0 and variance 1.
For $v\geq 0$ with $w=\max\{v,1\}$ the inequality
\begin{eqnarray*}
P(Z\geq v)\leq \frac{1}{2w^2}
\end{eqnarray*}
holds, with equality if
$$Z=\left\{\begin{array}{rcl}
-w & & \frac 1{2w^2}\\
\\
 0 & \mbox{with probability} & 1-\frac 1{w^2}\\
\\
w & & \frac 1{2w^2}\,.\\
\end{array}\right.$$
\end{theorem}\medskip

In view of $2w^2=2(\max\{v,1\})^2 \geq 1+v^2$ this inequality improves the bound from Cantelli's inequality from Theorem \ref{Cantelli}, as it
should.
For symmetric random variables we obtain the following bound for asymmetric intervals.

\begin{theorem}\label{symmetricgeneral}
Let the standardized random variable $Z$ be symmetric.
Consider any positive $u$ and $v$ with $v \leq u$ and discern four cases. \\
For $0 < v \leq u \leq 1$ the inequality
\begin{equation}\label{symmetricinequality11}
 P\left( Z \leq -u\ \  \mbox{\rm or}\ \ Z \geq v \right) \leq 1
\end{equation}
holds with equality if $Z$ puts mass $1/2$ at both $1$ and $-1$. \\
For $0 < v \leq u \leq {\sqrt2}\, v, 1 \leq u,$ the inequality
\begin{equation}\label{symmetricinequality2}
P\left( Z \leq -u\ \  \mbox{\rm or}\ \  Z \geq v \right) \leq \frac 1{u^2}
\end{equation}
is valid with equality if
\begin{eqnarray*}
Z = \left\{ \begin{array}{rcl}
            -u & & \frac 1{2 u^2}\\
            \\
            0 & \mbox{with probability} & 1-\frac 1{u^2}\\
            \\
            u & & \frac 1{2 u^2}\\
        \end{array}\right.
\end{eqnarray*}
holds.
For $0 < v \leq 1 < u, {\sqrt2}\, v \leq u,$ the inequality
\begin{equation}\label{symmetricinequality3}
P\left( Z \leq -u\ \  \mbox{\rm or}\ \  Z\geq v \right) \leq \frac 12 + \frac{1-v^2}{2(u^2-v^2)}
\end{equation}
is valid with equality if
\begin{eqnarray*}
Z = \left\{\begin{array}{rcl}
            -u,\, u & & \frac {1 - v^2}{2(u^2 - v^2)} \\
            & \mbox{with probability} & \\
            -v,\, v & & \frac {u^2 - 1}{2(u^2 - v^2)}
        \end{array}\right.
\end{eqnarray*}
holds.
For $1 \leq v \leq u, {\sqrt2}\, v \leq u,$ the inequality
\begin{equation}\label{symmetricinequality4}
P\left( Z \leq -u\ \  \mbox{\rm or}\ \  Z\geq v \right) \leq \frac 1{2 v^2}
\end{equation}
is valid with equality if
\begin{eqnarray*}
Z = \left\{\begin{array}{rcl}
            -v & & \frac 1{2 v^2}\\
            \\
            0 & \mbox{with probability} & 1-\frac 1{v^2}\\
            \\
            v & & \frac 1{2 v^2}\\
      \end{array}\right.
\end{eqnarray*}
holds.
\end{theorem}
Note that any choice of $(u,v)$ with $0<v \leq u$ belongs to at least one of the four cases in this Theorem.\\
\noindent
{\bf Proof}\\
To prove these inequalities we determine the supremum of the left hand side of (\ref{symmetricinequality11}) over all symmetric random
variables $Z$ with mean 0 and variance {\em at most} 1.
Let $Z$ be such a random variable and define the symmetric random variable $Y$ by
\begin{eqnarray*}
Y= -u{\bf 1}_{[Z\leq -u]} -v{\bf 1}_{[-u<Z\leq -v]} + v{\bf 1}_{[v\leq Z < u]} + u{\bf 1}_{[u \leq Z]}.
\end{eqnarray*}
Note that $Y$ is a discrete, symmetric random variable with probability mass at ${\cal V} = \{-u, -v, 0, v, u\}$ only and with
${\rm var}(Y) = E(Y^2) \leq E(Z^2) = {\rm var}(Z) \leq 1$.
Furthermore,
\begin{equation}\label{equality}
P\left( Y \leq -u\ \  \mbox{\rm or}\ \  Y \geq v \right) = P\left( Z \leq -u\ \  \mbox{\rm or}\ \ Z \geq v \right)
\end{equation}
holds, and we may conclude that the supremum of (\ref{equality}) over $Z$ is attained by a symmetric discrete random variable $Y$
taking its values at $\cal V$ and with $E(Y^2) \leq 1$.

We introduce
\begin{eqnarray*}
p = P(Y= -u) = P(Y=u), \quad q = P(Y=-v) = P(Y=v),
\end{eqnarray*}
and note that the supremum of (\ref{equality}) equals the maximum of
\begin{eqnarray*}
P\left( Y \leq -u\ \  \mbox{\rm or}\ \  Y \geq v \right) = 2p+q \nonumber
\end{eqnarray*}
over the convex polygon
\begin{eqnarray*}
{\cal Q} = \{(p,q) \mid p \geq 0 \ , \ q \geq 0 \ , \ p+q \leq \tfrac 12 \ , \ u^2 p + v^2q \leq \tfrac 12 \}. \nonumber
\end{eqnarray*}
In this linear programming problem the maximum is attained at one of the vertices of polygon $\cal Q$.
We discern three cases.
\begin{itemize}
\item[\bf A.]$\ 0 < v \leq u \leq 1$ \\
Here, $\cal Q$ reduces to the triangle $$ {\cal Q} = \{(p,q) \mid p \geq 0 \ , \ q \geq 0 \ , \ p+q \leq \frac12 \},$$
the maximum 1 of the map $(p,q) \mapsto 2p+q$ on $\cal Q$ is attained at $(p,q) = (\tfrac 12,0)$, and we get inequality (\ref{symmetricinequality11}).

\item[\bf B.]$\ 0 < v \leq 1 < u$ \\
In this case the polygon $\cal Q$ is a quadrangle with vertices
\begin{equation}
(0, 0), \left(\frac 1{2u^2}, 0\right), \left(0, \tfrac 12 \right), \left(\frac{1-v^2}{2(u^2-v^2)}, \frac{u^2-1}{2(u^2-v^2)} \right).
\end{equation}
The corresponding values of the function $(p,q) \mapsto 2p+q$ are
\begin{eqnarray*}
0, \ \ \frac 1{u^2}, \ \ \frac 12, \ \ \frac 12 + \frac{1-v^2}{2(u^2-v^2)}.
\end{eqnarray*}
Computation shows that the fourth value is larger than the second value and hence largest, iff ${\sqrt 2}\, v \leq u$ holds.
Note that this yields inequality (\ref{symmetricinequality3}) and inequality (\ref{symmetricinequality2}) under the additional restriction
$v \leq 1$.

\item[\bf C.]$\ 1 \leq v \leq u$ \\
Polygon $\cal Q$ reduces to a triangle here with vertices
\begin{eqnarray*}
(0, 0), \left(0, \frac 1{2v^2}\right), \left( \frac 1{2u^2}, 0 \right).
\end{eqnarray*}
The corresponding values of the function $(p,q) \mapsto 2p+q$ are
\begin{eqnarray*}
0, \ \ \frac 1{2v^2}, \ \ \frac 1{u^2}.
\end{eqnarray*}
Computation shows that this implies inequality (\ref{symmetricinequality4}) and inequality (\ref{symmetricinequality2})
under the additional restriction $v \geq 1$.
\end{itemize}
Straightforward computation shows that equalities are attained by the random variables mentioned in the Theorem.
\hfill$\Box$\\

\section{Concave distribution functions}\label{concave}
%\label{concave}

The upper bounds from the preceding sections to the probability in $(\ref{rint1})$ are rather large.
It is to be expected that restriction of the class of completely unknown distributions and the class of symmetric distributions to smaller
classes of distributions will yield smaller upper bounds.
In the next three sections we will obtain sharp upper bounds for $(\ref{rint1})$ over the class of unimodal distributions, the class of unimodal distributions with mean and mode coinciding, and the class of symmetric unimodal distributions, respectively.
This unimodality assumption is not unrealistic, as it is a very natural assumption in several practical applications, like statistical process control.

A distribution is unimodal with mode at $M$ if its corresponding distribution function is convex on $(-\infty,M)$ and concave on $[M,\infty).$ Consequently, a unimodal distribution has at most one atom, which may occur only at the mode $M$.
If a unimodal distribution is uniform on its support with an atom at one of its boundary points, we will call it a one-sided boundary-inflated
uniform distribution; cf. \cite{Klaassen}.
We shall repeatedly use a representation theorem for unimodal distributions of Khintchine, Lemma \ref{Khintchinelemma}, \cite{Khintchine}.
It characterizes unimodal distributions as a mixture of uniform distributions.
The inequalities we will derive attain equality for mixtures of at most three uniforms, where often one of these uniforms is degenerate, i.e., a point mass.
Unimodal distributions with their mode at $M=0$ and all their mass on the nonnegative half line $[0, \infty)$ have a distribution function that is concave on $[0,\infty)$ and vanishes on $(-\infty,0).$
They have a nonincreasing density on $(0,\infty).$
This special class of distributions is considered in the present section.
For this class a one-sided version of the Gau\ss \, inequality holds.
The Gau\ss \, inequality itself is an immediate consequence of it and will be presented also.

\begin{theorem}\label{onesidedGauss}{\bf One-sided Gau\ss \, inequality}\\
Let the random variable $Y$ have second moment $E(Y^2) =1$ and let its distribution function be concave on $[0,\infty)$ and $0$ on $(-\infty,0).$ For all nonnegative $v$ the inequality
\begin{eqnarray*}
P(Y \geq v) \leq {\left\{
\begin{array}{rcl}
    1- \frac v{\sqrt 3}, &  & 0 \leq v \leq \frac 2{\sqrt 3},\\
    \\
    \frac 49 \, \frac 1{v^2},& & \frac 2{\sqrt 3} \leq v, \\
\end{array}\right.}
\end{eqnarray*}
is valid.
For $0 \leq v \leq 2/{\sqrt 3}$ equality holds iff $Y$ has a uniform distribution on $[0,\sqrt 3).$
For $2/{\sqrt 3} \leq v$ equality holds iff $Y$ has a one-sided boundary-inflated
uniform distribution on $[0, 3v/2)$ with mass $1 -4/(3v^2)$ at $0.$
\end{theorem}
{\bf Proof}\\
By Khintchine's representation from Lemma \ref{Khintchinelemma} there exist a probability $p_0$ and a distribution function $F$ on the positive half line, such that $P(Y=0)=p_0$ holds and the density of $Y$ at $y$ on the positive half line equals
$(1-p_0)\int_0^\infty c^{-1}{\bf 1}_{(0,c)}(y)\, dF(c).$
By Fubini it follows that
\begin{equation}\label{concave2}
1=E(Y^2) = (1-p_0)\int_0^\infty \frac 13 c^2\, dF(c)
\end{equation}
holds and that for positive $v$
\begin{eqnarray}\label{concave3}
\lefteqn{P(Y \geq v) = (1-p_0) \int_0^\infty \int_v^\infty \frac 1c {\bf 1}_{(0,c)}(y)\, dy\, dF(c) \nonumber}\\
&& =(1-p_0) \int_v^\infty \left(1 - \frac vc \right) dF(c)
\end{eqnarray}
is valid.
Without loss of generality we assume that $F$ puts positive mass on $(v, \infty).$
Let us write
\begin{eqnarray*}
c_v = \int_v^\infty c\, dF(c)/(1-F(v)) > v.
\end{eqnarray*}
As the map $c \mapsto (1-v/c)$ is strictly concave on $[v,\infty)$, (\ref{concave3}) implies by Jensen's inequality
\begin{equation}\label{concave5}
P(Y \geq v) \leq (1-p_0) (1-F(v))\left(1-\frac v{c_v} \right) \leq (1-p_0)\left(1-\frac v{c_v} \right)
\end{equation}
with equalities iff $F$ is degenerate at $c_v.$ This means that we may restrict attention to those $Y$ with mass $p_0$ at 0 for which $F$ is degenerate at some $c_v > v.$ For such $Y$ equation (\ref{concave2}) implies
\begin{equation}\label{concave6}
1-p_0 = 3c_v^{-2},
\end{equation}
which together with (\ref{concave5}) yields
\begin{equation}\label{concave7}
P(Y \geq v) \leq \frac 3{c_v^2} \left(1-\frac v{c_v} \right).
\end{equation}
As this function of $c_v$ attains its maximum at $c_v = 3v/2 > v,$ we obtain
\begin{eqnarray*}
P(Y \geq v) \leq \frac 4{9v^2}
\end{eqnarray*}
with equality iff $Y$ has the one-sided boundary-inflated uniform distribution as described in the Theorem.
However, for $c_v = 3v/2$ equation (\ref{concave6}) becomes $1-p_0 = 4/(3v^2),$ which for $v < 2/\sqrt 3$ leads to an impossible, negative value of $p_0.$
This means that for $v < 2/\sqrt 3$ the mass at 0 vanishes and that (\ref{concave5}), (\ref{concave6}), and (\ref{concave7}) hold with $c_v = \sqrt 3.$
\hfill$\Box$\\

About two centuries ago Johann Carl Friedrich Gau\ss \, presented and proved a sharp upper bound on the probability $P(\lvert X \rvert  \geq v)$ for unimodal random variables $X$ with mode at 0 and finite second moment in Sections 9 and 10 of \cite{Gauss}; for a translation from Latin into English see \cite{Stewart}.
His result precedes the famous Bienaym\'e-Chebyshev inequality (\ref{BCi}) by three decades.
The Gau\ss \, inequality for large values of $v$ has been given in (\ref{Gauss0}).
The complete inequality is the following.

\begin{theorem}\label{originalGauss}{\bf Original Gau\ss \, inequality}\\
Let the random variable $Y$ have a unimodal distribution with mode at 0 and second moment $E(Y^2) =1$.
For all nonnegative $v$ the inequality
\begin{eqnarray}\label{originalGauss1}
P(\lvert Y \rvert   \geq v) \leq {\left\{
\begin{array}{rcl}
    1- \frac v{\sqrt 3}, &  & 0 \leq v \leq \frac 2{\sqrt 3},\\
    \\
    \frac 49 \, \frac 1{v^2},& & \frac 2{\sqrt 3} \leq v, \\
\end{array}\right.}
\end{eqnarray}
is valid. For $0 \leq v \leq 2/{\sqrt 3}$ equality holds if $Y$ has a uniform distribution on $(-\sqrt{3},\sqrt{3}).$
For $2/{\sqrt 3} \leq v$ equality holds if $Y$ has mass $1 -4/(3v^2)$ at $0$
and the rest of its mass uniformly distributed on $(-3v/2, 3v/2)$.
\end{theorem}
{\bf Proof}\\
As $\lvert Y \rvert  $ has a concave distribution function on $[0,\infty)$, the one-sided Gau\ss \, inequality proves the Theorem.
\hfill$\Box$\\

Our proof of the Gau\ss \, inequality via the Khintchine representation and Jensen's inequality differs from the three proofs as presented by \cite{Pukelsheim}.
Observe that the bound in (\ref{originalGauss1}) can also be described as the minimum of the two functions in there.
Also note that Gau\ss \, considered only densities and hence could not prove the second bound in (\ref{originalGauss1}) to be sharp.

\section{Unimodal distributions}\label{unimodal}

In the preceding section on concave distributions we have already defined the related class of unimodal distributions, which we will study in the remaining sections. The factor 4/9 from the one-sided Gau\ss \, inequality of Theorem \ref{onesidedGauss} will play a role in all these sections.
For the proof of our extension of the Cantelli inequality from Theorem \ref{Cantelli} to unimodal distributions we shall use the following powerful result for unimodal distributions, which also shows the factor 4/9.

\begin{theorem}\label{VPtheorem}{\bf Vysochanski${\rm {\bf \breve{i}}}$ and Petunin inequality}\\
Any unimodal random variable $W$ with finite second moment satisfies
\begin{eqnarray}\label{VP1}
P(\lvert W \rvert   \geq w) \leq {\left\{
\begin{array}{rcl}
    1 & & 0 \leq w \leq \sqrt{E(W^2)}, \\
    \\
    \frac 43\, \frac{E(W^2)}{w^2} - \frac 13 & {\rm for} & \sqrt{E(W^2)} \leq w \leq \sqrt{8E(W^2)/3}, \\
    \\
    \frac 49\, \frac{E(W^2)}{w^2} & & \sqrt{8E(W^2)/3} \leq w.\\
\end{array}\right.}
\end{eqnarray}
\end{theorem}

\noindent
{\bf Proof}\\
The proof of \cite{Vysochanskii1} and \cite{Vysochanskii2} has been smoothed by \cite{Pukelsheim} and invokes Gau\ss's inequality presented in Theorem \ref{originalGauss}.
\hfill$\Box$\\

Actually, for $\sqrt{8/3} \leq v$ inequality (\ref{VP1}) implies the Gau\ss \, inequality (\ref{originalGauss1}).
Observe that the bound in (\ref{VP1}) can be described as the minimum of the three expressions at its right hand side.

Here is our analogue of Cantelli's inequality from Theorem \ref{Cantelli}.

\begin{theorem}\label{unimodalCantelli}
Let the distribution of the standardized random variable $Z$ be unimodal. For any
$v\geq 0$ the inequality
\begin{equation}\label{unimodalCantellia}
P(Z\geq v)\leq {\left\{
\begin{array}{rcl}
    \frac {3 - v^2}{3(1 + v^2)}, & & 0 \leq v \leq \sqrt{\frac 53},\\
    \\
    \frac 49 \, \frac 1{1 + v^2},& & \sqrt{\frac 53} \leq v, \\
\end{array}\right.}
\end{equation}
holds, with equality for $0 < v \leq \sqrt{5/3}$ if $Z$ has mass $(3 - v^2)/(3(1+v^2))$ at $v$ and the rest of its mass, $4 v^2/(3(1+v^2)),$ uniformly distributed on the interval $[-(3+v^2)/(2v),v]$, and with equality for $\sqrt{5/3} \leq v$ if $Z$ has mass $(3v^2-1)/(3(1+v^2))$ at
$-1/v$ and the rest of its mass, $4/(3(1+v^2)),$ uniformly distributed on the interval $[-1/v, (1+3v^2)/(2v)].$
\end{theorem}

\noindent
{\bf Proof}\\
Applying Theorem \ref{VPtheorem} with $W=Z + 1/v$ and $w=v + 1/v$ we obtain (\ref{unimodalCantellia}) after some computation.
Additional computation shows that the random variables mentioned in the theorem attain equality.
\hfill$\Box$\\

Comparing this inequality (\ref{unimodalCantellia}) to Cantelli's inequality from Theorem \ref{Cantelli} we note the extra factor 4/9 for larger values of $v$; see also Table \ref{table1}.
Furthermore, note that the bound in (\ref{unimodalCantellia}) can be viewed as the minimum of the two functions in there, and that these functions intersect at $v=\sqrt{5/3}$.

Next we turn to the general case of asymmetric intervals around 0.
The Vysochanski${\rm \breve{i}}$ and Petunin inequality from Theorem \ref{VPtheorem} implies the following result.

\begin{theorem}\label{generalGaussa}
For $v \geq \sqrt{5/3},\ \max \{v,(11v -4 \sqrt{6v^2-10})/5 \} \leq u \leq v+2/v$, and any standardized unimodal random variable $Z$, the inequality
\begin{equation}\label{boundgeneralGaussa}
P(Z\leq -u\ \  \mbox{\rm or}\ \  Z \geq v) \leq \frac 49\ \frac{4 + (u-v)^2}{(u+v)^2}
\end{equation}
holds with equality if $Z=(v-u)/2 + UY$, $U$ and $Y$ independent random variables, $U$ uniform on the unit interval, and $Y$ the generalized Bernoulli random variable
\begin{eqnarray}\label{defY3}
Y = \left\{
\begin{array}{rcl}
-\frac 34 (u+v) &                             & \frac 43 \frac{2 + v(v-u)}{(u+v)^2} \\
0               & \mbox{with probability} & 1 - \frac 43 \frac{4 + (v-u)^2}{(u+v)^2} \\
\frac 34 (u+v)  &                             & \frac 43 \frac{2 + u(u-v)}{(u+v)^2}.
\end{array}\right.
\end{eqnarray}
\end{theorem}

\noindent
{\bf Proof}\\
As in the proof of Theorem \ref{general} we note
\begin{eqnarray*}
P(Z\leq -u\ \ \mbox{or}\ \  Z \geq v) = P\left(\lvert 2Z - (v-u) \rvert   \geq u+v \right).
\end{eqnarray*}
Applying the third inequality of (\ref{VP1}) from Theorem \ref{VPtheorem} we obtain (\ref{boundgeneralGaussa}).
Computation shows that the random variable $Y$ and hence $Z = (v-u)/2 + UY$ are well defined under the conditions on $u$ and $v$, and that this $Z$ attains the bound.
\hfill$\Box$\\

An immediate consequence of this Theorem is the following one, which is the main content of Theorem 2 of \cite{Vysochanskii1}.

\begin{corollary}\label{absoluteGauss}
For $\sqrt{8/3} \leq v$ and any standardized unimodal random variable $Z$, the inequality
\begin{eqnarray*}
P(\lvert Z \rvert   \geq v) \leq \frac 4{9v^2}
\end{eqnarray*}
holds with equality for $Z=UY$, $U$ and $Y$ independent random variables, $U$ uniform on the unit interval, and $Y$ the generalized Bernoulli random variable
\begin{eqnarray*}
Y = \left\{
\begin{array}{rcl}
-\frac 32 v &                             & \frac 2{3v^2} \\
0           & \mbox{with probability} & 1 - \frac 4{3v^2} \\
\frac 32 v  &                             & \frac 2{3v^2}.
\end{array}\right.
\end{eqnarray*}
\end{corollary}

Instead of applying the Vysochanski${\rm \breve{i}}$ and Petunin inequality from Theorem \ref{VPtheorem} we could choose the approach via Khintchine's characterization of unimodal distributions and Jensen's inequality as in Section \ref{concave}.
An admittedly laborious proof along these lines of Cantelli's inequality for unimodal distributions as given in Theorem \ref{unimodalCantelli}, is presented in Subsection \ref{proofunimodalCantelli} of the Appendix.
However, this Khintchine-Jensen approach yields a partially improved version of Theorem \ref{generalGaussa} too, namely
\begin{theorem}\label{generalGauss}
Assume $\sqrt{3} \leq v \leq u$. For any standardized unimodal random variable $Z$, the inequality
\begin{eqnarray}\label{boundgeneralGauss}
P(Z\leq -u\ \  \mbox{\rm or}\ \  Z \geq v) \leq \left\{
\begin{array}{rcl}
\frac 49\ \frac{4 + (u-v)^2}{(u+v)^2} &                & v \leq u \leq v+ \frac 2v \\
                                      & \mbox{\rm for} &  \\
\frac 49\ \frac 1{1+v^2}              &                & v + \frac 2v \leq u
\end{array}\right.
\end{eqnarray}
holds.
For $v \leq u \leq v + \frac 2v$ equality is attained if $Z$ is defined as in (\ref{defY3}).
For $v + \frac 2v \leq u$ equality is attained if $Z$ has mass $(3v^2-1)/(3(1+v^2))$ at
$-1/v$ and the rest of its mass, $4/(3(1+v^2)),$ uniformly distributed on the interval $[-1/v, (1+3v^2)/(2v)]$, like in Theorem \ref{unimodalCantelli}.
\end{theorem}

The proof of this Theorem is given in Subsection \ref{proofgeneralGauss} of the Appendix.

\section{Unimodal distributions with coinciding mode and mean}\label{unimodalcoinciding}

When we restrict the class of distributions further to the class of unimodal distributions with coinciding mode and mean, then for the one-sided probability we see that the factor $4/9$ does not play such a role anymore as in Theorem \ref{unimodalCantelli}, the analogue of Cantelli's inequality, Theorem \ref{Cantelli}.

\begin{theorem}\label{onesidedGauss2}
For any standardized unimodal random variable $Z$ with mode at 0, the inequality
\begin{eqnarray}\label{boundonesidedGauss}
\lefteqn{ P(Z \geq v)\leq \frac {2(x-1)}{v^2x^2 + 2x +1}, }\\
&& x= \frac 12 \left( w + 1 + \frac 1w \right), \quad w= \left( \sqrt{3+v^2} + \sqrt{3} \right)^{2/3} v^{-2/3}, \nonumber
\end{eqnarray}
holds with equality if $Z=UY$ holds with $U$ and $Y$ independent random variables, $U$ uniform on the unit interval, and $Y$ the Bernoulli variable
\begin{eqnarray}\label{defY0}
Y = \left\{
\begin{array}{rcl}
\frac{-3}{vx} & & \frac{v^2 x^2}{3 + v^2 x^2}\\
   &\mbox{with probability}& \\
vx & & \frac 3{3 + v^2 x^2}.
\end{array}\right.
\end{eqnarray}
\end{theorem}

\noindent
{\bf Proof}\\
We apply the preparatory part of the alternative proof of Theorem \ref{unimodalCantelli} in the Appendix with $M=0$.
We see that the supremum at the left hand side of (\ref{Ppsi}) with $M=0$ is given by (\ref{supg}) with $M=0$ and equals
\begin{equation}\label{supg2}
\sup_{y \geq 0} \frac y{1+y}\ \frac 3{3 + v^2 (1+y)^2} = \sup_{x \geq 1} \frac {x-1}x \ \frac 3{3 + v^2 x^2}.
\end{equation}
Straightforward computation shows that the derivative with respect to $x$ of the function in (\ref{supg2}) is nonnegative if and only if
\begin{equation}\label{derivativepositive1}
2v^2 x^3 - 3v^2 x^2 - 3 \leq 0
\end{equation}
holds.
Observe that the function $x \mapsto 2x^3 -3x^2$ is increasing for $x \geq 1$ and negative for $1 \leq x < 3/2$.
By Vieta's method to tackle cubic equations we substitute $x=(w+1 + 1/w)/2 \geq 3/2, \ w > 0,$ and obtain equality in (\ref{derivativepositive1})
if and only if
\begin{equation*}
v^2 \left( w^3 \right)^2 - 2\left( 6 + v^2 \right) w^3 + v^2 = 0
\end{equation*}
holds.
The positive roots of this quadratic equation in $w^3$ yield $w_1 = (\sqrt{3+v^2} + \sqrt{3})^{2/3} v^{-2/3}$ and
$w_2 = (\sqrt{3+v^2} - \sqrt{3})^{2/3} v^{-2/3}$.
Note that $w_1 w_2 =1$ and hence $w_1 + 1/w_1 = w_2 + 1/w_2$ hold.
Consequently, the only real root $x$ of the cubic function in (\ref{derivativepositive1}) equals the one given in (\ref{boundonesidedGauss}).
Combining (\ref{supg2}) and (\ref{derivativepositive1}) (with equality) we arrive at the inequality in (\ref{boundonesidedGauss}). \\
Straightforward verification shows that $Z=UY$ with $Y$ as in (\ref{defY0}) attains this bound.
\hfill $\Box$\\

The sharp, restricted Gau\ss \, inequality for random variables with coinciding mean and mode is the same as the original one from Theorem \ref{originalGauss}, as the distributions that attain equality in (\ref{originalGauss1}) are symmetric and hence have coinciding mean and mode.

\begin{theorem}\label{Gauss}{\bf Restricted Gau\ss \, inequality}\\
For any standardized unimodal random variable $Z$ with mode at 0, the inequality
\begin{eqnarray}\label{Gaussineq}
P(\lvert Z \rvert   \geq v)\leq  {\left\{\begin{array}{rcl}
1- \frac v{\sqrt 3}, &  & 0\leq v < \frac 2{\sqrt 3},\\
\\
\frac 49 \, \frac 1 {v^2},& & v \geq \frac 2{\sqrt 3}, \\
\end{array}\right.}
\end{eqnarray}
holds with equality if the distribution of $Z$ is the mixture of a
uniform distribution on $[-((\frac{3}{2}v)\vee \sqrt 3)\ ,\
(\frac{3}{2}v)\vee \sqrt 3]$ and a distribution degenerate at $0$
such that the point mass at $0$ equals $[1-4/(3v^2)]\vee 0$.
\end{theorem}

We extend Gau\ss's inequality to asymmetric intervals as in (\ref{rint1}) as follows.
\begin{theorem}\label{twosidedGauss2}
For $\sqrt{3} \leq u \leq v \leq u + 2/u$ or $\sqrt{3} \leq v \leq u \leq v + 2/v$ and for any standardized unimodal random variable $Z$ with mode at 0, the inequality
\begin{equation}\label{boundtwosidedGauss}
P(Z\leq -u\ \  \mbox{\rm or}\ \  Z \geq v) \leq \frac{1024 v^{-2} + 27(\gamma -1)^2(\gamma +3)^2}{9(\gamma +3)^3(3\gamma +1)}
\end{equation}
holds with
\begin{equation}\label{gamma0}
\gamma = 2 \sqrt{1+\frac uv} \cos \left( \frac 13 \left[\pi -\arctan \left(\sqrt{\frac uv}\right)\right] \right) - 1.
\end{equation}
Equality is attained in (\ref{boundtwosidedGauss}) for $Z=UY$ with $U$ and $Y$ independent random variables, $U$ uniform on the unit interval, and $Y$ the generalized Bernoulli random variable
\begin{eqnarray}\label{defY3M0}
Y = \left\{
\begin{array}{rcl}
- \gamma c            &                             & 3 \frac{4c^{-2} - \gamma +1}{(\gamma +1)(\gamma +3)} \\
\frac 34 (1-\gamma)c & \mbox{with probability}  & 16 \frac{\gamma -3c^{-2}}{(\gamma +3)(3\gamma +1)} \\
c                     &                             & 3 \frac{4c^{-2} + \gamma(\gamma -1)}{(\gamma +1)(3\gamma +1)}
\end{array}\right.
\end{eqnarray}
with
\begin{equation}\label{valuec}
c = \frac 38 (\gamma +3) v.
\end{equation}
\end{theorem}

\noindent
\noindent {\bf Proof}\\
The proof of Theorem \ref{generalGauss}, given at the end of the Appendix, can be applied with $M=0$ all the way up to and including the value of the upper bound in (\ref{LA620a}) with $\gamma$ defined in (\ref{gamma0}) and satisfying
\begin{equation}\label{gammaM0}
\gamma^3 + 3\gamma^2 - 3\frac uv \gamma - \frac uv = 0.
\end{equation}
Note that the equation (\ref{gammaM0}) may be rewritten as
\begin{equation*}
(3\gamma +1) u = (\gamma^3 + 3 \gamma^2) v,
\end{equation*}
which implies
\begin{equation*}
(3\gamma +1)(u+v) = (\gamma +1)^3 v.
\end{equation*}
With the help of the last two equations $u$ may be eliminated from (\ref{LA620a}) resulting in the expression in (\ref{boundtwosidedGauss}).

In view of this the random variable $Y$ from (\ref{defY3M0}) follows straightforwardly from (\ref{LA618a}) and (\ref{pandq}) with $M=0$.
\hfill$\Box$\\

\begin{remark}\label{u=v}
For $u=v$ the value of $\gamma$ from (\ref{gamma0}) becomes 1 and the upper bound in (\ref{boundtwosidedGauss}) takes on the value $4/(9v^2)$, which is in line with Theorem \ref{Gauss}.
\end{remark}

\section{Symmetric Unimodal distributions}\label{symmetricunimodal}

Under the extra assumption of unimodality, Theorem \ref{symmetricChebyshev} for symmetric distributions may be sharpened too.
Again we will encounter the extra factor $4/9$.
The resulting inequality for symmetric unimodal distributions has been obtained by \cite{Camp}, \cite{Meidell}, and \cite{Shewhart}.
A different proof is given by \cite{Theil}.
Still a different proof is given in Lemma 2 of \cite{Clarkson}.
However, our proof is shorter and simpler than theirs.
\begin{theorem}\label{Gausslike}
If the distribution of the standardized random variable $Z$ is symmetric and unimodal, then
\begin{eqnarray*}
P(Z \geq v)\leq  {\left\{\begin{array}{rcl}
\frac 12 \left( 1- \frac v{\sqrt 3} \right), &  & 0 < v < \frac 2{\sqrt 3},\\
\\
\frac 49 \, \frac 1 {2v^2},& & v \geq \frac 2{\sqrt 3}, \\
\end{array}\right.}
\end{eqnarray*}
holds. Equality is attained by the mixture of a uniform
distribution on\\
 $[-((\frac{3}{2}v)\vee \sqrt 3)\ ,\ (\frac{3}{2}v)\vee \sqrt 3]$ and a distribution degenerate at $0$
such that the point mass at $0$ equals $[1-4/(3v^2)]\vee 0$.
\end{theorem}
{\bf Proof}\\
As $Z$ is symmetric, $P(Z \geq v) = P(\lvert Z \rvert   \geq v)/2$ holds and hence Theorem \ref{Gauss} yields a proof.
\hfill$\Box$\\

Actually, this Theorem is equivalent to Theorem \ref{Gauss}.
Indeed, let $Z$ be standardized with mode at 0 and let $B$ be an independent Bernoulli random variable with $P(B=-1) =P(B=1) = 1/2.$
As $BZ$ is symmetric and hence $P(\lvert Z \rvert   \geq v) = P(\lvert BZ \rvert   \geq v) = 2P(BZ \geq v)$ holds, Theorem \ref{Gausslike} implies Theorem \ref{Gauss}.
For the class of symmetric unimodal distributions these Theorems also imply that inequality (\ref{Gaussineq}) holds and is sharp.

Next we consider the case of asymmetric intervals around 0.

\begin{theorem}\label{symmetricunimodaltheorem}
Assume the distribution of the standardized random variable $Z$ is symmetric and unimodal and consider $\sqrt 3 \leq v \leq u.$
For $u \leq (2 \sqrt 2 -1)v$ the inequality
\begin{equation}\label{su1}
P(Z\leq -u\ \  \mbox{\rm or}\ \  Z \geq v) \leq \frac 49 \ \frac 4{(u+v)^2}
\end{equation}
holds with equality if $Z$ is uniform on $[- 3(u+v)/4\,,\, 3(u+v)/4]$ with probability $16/(3(u+v)^2)$
and has a point mass at 0 with probability $1-16/(3(u+v)^2).$

For $(2 \sqrt 2 -1)v \leq u$ the inequality
\begin{equation}\label{su2}
P(Z\leq -u\ \  \mbox{\rm or}\ \  Z \geq v) \leq \frac 49 \ \frac 1{2v^2}
\end{equation}
holds with equality if $Z$ is uniform on $[- 3v/2\,,\, 3v/2]$ with probability $4/(3v^2)$ and has a point mass at 0 with probability $1-4/(3v^2).$
\end{theorem}

\begin{remark}\label{Semenikhin}
The Theorem in Section 3 of \cite{Semenikhin} presents a complete version for all positive $u$ and $v$ of our Theorem \ref{symmetricunimodaltheorem} with his $m=(v-u)/2$ and $h=(u+v)/2$.
\cite{Semenikhin} uses another approach for the proof than we do, although there are similarities between (\ref{eq1}) below and his expression in the displayed formula above his (A.3).
\end{remark}

\noindent
{\bf Proof}\\
By Lemma \ref{Khintchinelemma} we have $Z=UY$ with $Y$ symmetric around 0 because of the symmetry of $Z$ around 0.
Along the lines of Lemma \ref{unimodalCantelliPsi} we obtain by Jensen's inequality
\begin{eqnarray*}
\lefteqn{P(Z\leq -u\ \  \mbox{\rm or}\ \  Z \geq v) = P(Z\geq u) + P(Z \geq v) \nonumber} \\
& =& E\left(1- \frac vY \mid v \leq Y < u \right)\,P(v \leq Y < u) \nonumber \\
& & + E\left(2- \frac {u+v}Y \mid Y \geq u \right)\,P(Y\geq u)\\
& \leq &  \left[1- \frac v{\sqrt{E(Y^2 \mid v \leq Y < u)}}\right] \, P(v \leq Y < u) \nonumber \\
& & + \left[2- \frac {u+v}{\sqrt{E(Y^2 \mid Y\geq u)}}\right]\,P(Y\geq u). \nonumber
\end{eqnarray*}
With the notation $a = \sqrt{E(Y^2 \mid v \leq Y < u)},\, b = \sqrt{E(Y^2 \mid Y\geq u)}$, \\
$p = P(v \leq Y < u),$ and $q = P(Y\geq u)$ this implies
\begin{eqnarray}\label{eq1}
\lefteqn{ P(Z\leq -u\ \  \mbox{\rm or}\ \  Z \geq v)}\\
&& = \sup \left\{ \left(1-\frac va \right)\,p + \left(2-\frac {u+v}b \right)\,q \ :\ v \leq a \leq u \leq b, \right. \nonumber \\
&& \qquad \qquad \qquad \qquad \qquad \left. p \geq 0,\ q \geq 0,\ p+q \leq \tfrac 12,\ 2a^2p + 2b^2q \leq 3 \right\}. \nonumber
\end{eqnarray}
By increasing $a$ or $b$ if necessary, we see that this supremum is attained at \\
$2a^2p + 2b^2q =3.$

Fix $a$ and $b$ with $v \leq a \leq u \leq b$ and write $\alpha = 1- v/a \geq 0$ and $\beta= 2- (u+v)/b \geq 0.$
Consider
\begin{eqnarray*}
S(\alpha, \beta) = \sup \{ \alpha p + \beta q\,:\,p \geq 0,\ q \geq 0,\ p+q \leq \tfrac 12,\ 2a^2p + 2b^2q = 3 \}.
\end{eqnarray*}
Now $b \geq a \geq v \geq \sqrt 3$ and hence $p+q \leq (a^2p + b^2q)/3 = 1/2$ hold.
Consequently, we have
\begin{equation}\label{eq2}
S(\alpha, \beta) = \sup \left\{ \alpha p + \beta \, \frac {3-2a^2p}{2b^2} \,:\, 0 \leq p \leq \frac 3{2a^2} \right\}
= \frac 32 \max \left\{ \frac \alpha {a^2}, \frac \beta {b^2} \right\}.
\end{equation}
Studying the stationary points of the functions $a \mapsto (1- v/a)/a^2$ and \\
$b \mapsto (2- (u+v)/b)/b^2$ we see with the help of (\ref{eq2}) that the supremum in (\ref{eq1}) equals
\begin{eqnarray*}
\lefteqn{\frac 32 \max\left\{ \sup_{v \leq a \leq u}\left[\frac 1{a^2} \left(1 - \frac va \right) \right],
 \sup_{b \geq u} \left[\frac 1{b^2} \left(2 - \frac {u+v}b \right)\right] \right\} } \\
&& \qquad \qquad \qquad = \max \left\{ \frac 2{9v^2}{\bf 1}_{[3v \leq 2u]}\,,\ \frac 3{2u^2} \left(1-\frac uv \right),\
\frac {16}{9(u+v)^2} {\bf 1}_{[u \leq 3v]} \right\}.
\end{eqnarray*}
If $3v < u$ holds, then we have $$\frac 3{2u^2} \left(1 - \frac vu \right) \leq \frac 3{2u^2} < \frac 1{6v^2} < \frac 2{9v^2}.$$
If $3v \geq u \geq v$ holds, then $$\frac 3{2u^2} \left(1 - \frac vu \right) \leq \frac {16}{9(u+v)^2}$$ holds
iff $$\frac {27}{32} \left(1 - \frac vu \right) \left(1 + \frac vu \right)^2 \leq 1$$ holds.
The last inequality is valid in view of $$\sup_{1/3 \leq x \leq 1} \frac {27}{32} (1-x)(1+x)^2 = 1.$$
We conclude that the supremum in (\ref{eq1}) is bounded by $$\frac 49 \max \left\{ \frac 1{2v^2}\,,\ \frac 4{(u+v)^2} \right\}.$$
Straightforward computation shows that equality holds in (\ref{su1}) and (\ref{su2}) for the indicated random variables $Z.$
\hfill$\Box$\\

\appendix

\section{Appendix}

In this Appendix we prove a lemma with Khintchine's representation theorem, \cite{Khintchine}, and other lemmata that we need.
Furthermore, an alternative proof of Theorem \ref{unimodalCantelli} and the proof of Theorem \ref{generalGauss} are given.

\begin{lemma}{\bf Khintchine representation}\label{Khintchinelemma}\\
If $Z$ has a unimodal distribution, then there exist a constant
$M$ and independent random variables $U$ and $Y$ with $U$
uniformly distributed on the unit interval, such that $Z=M + UY$
holds. If $Z$ is symmetric around $M,$ then $Y$ is symmetric
around 0.
\end{lemma}
{\bf Proof}\\
Let $Z$ have its mode and possibly a point mass at $M$.
Theorem V.9 of \cite{Feller} and Theorem 1.3 of \cite{Dharmadhikari} yield the characterization $Z=M + UY.$
As an alternative proof, cf. page 8 of \cite{Dharmadhikari}, let the conditional distribution of $Z-M$ given $Z>M$ have
density $f$ and distribution function $F$ on
$(0,\infty).$ Since $f$ is nonincreasing, we may write
$$f(x)=-\int_0^\infty \frac 1y {\bf 1}_{[0, y]}(x)\, y\,df(y) =
\int_0^\infty \frac 1y {\bf 1}_{[0, y]}(x)\,d[F(y)-yf(y)]$$ with
$$F(y)-yf(y) = \int_0^y [f(z)-f(y)]\,dz$$ a distribution function
on $[0, \infty).$ It follows that $f$ is the density of
$UY$ given $Y$ positive. With a similar argument for negative
values of $Z-M$ and $Y$, we obtain Khintchine's characterization.
From this construction it follows that $Y$ is symmetric around 0,
if $Z$ is symmetric around $M$.
\hfill$\Box$\\

For the proofs of Theorems \ref{unimodalCantelli}, \ref{generalGauss}, \ref{onesidedGauss2}, \ref{twosidedGauss2}, and \ref{symmetricunimodaltheorem} we need also the following result.
\begin{lemma}\label{unimodalCantelliPsi}
Let ${\cal Z}_M$ be the class of random variables $Z$ that have a unimodal distribution with mean zero, unit
variance, and mode at $M.$ Let ${\cal Y}_M$ be the class of random variables $Y$ with mean $-2M$ and variance $3-M^2.$
With $-u \leq M \leq v$ the following holds true
\begin{eqnarray*}
\sup_{Z \in {\cal Z}_M}\ P(Z\leq -u\ \  \mbox{\rm or}\ \  Z \geq v) =
\sup_{Y \in{\cal Y}_M}\ E \Psi_M(Y)
\end{eqnarray*}
with the function $\Psi_M$ given by
\begin{eqnarray*}
\Psi_M(y)=\left[1+\frac{u+M}y \right] {\bf 1}_{[y<-u-M]}+ \left[
1-\frac{v-M}y \right] {\bf 1}_{[y>v-M]}.
\end{eqnarray*}
\end{lemma}
{\bf Proof}\\
By Lemma \ref{Khintchinelemma} we may write $Z=M+UY$ with $U$ and $Y$
independent, which implies $0=EZ=M+ \tfrac 12 EY$ and $1={\rm var}
Z = {\rm var}(UY)= E(U^2Y^2) - (E(UY))^2 = \tfrac 13 E(Y^2)- M^2.$
These equations yield $EY=-2M,\, E(Y^2) = 3(1+M^2),$ and hence ${\rm var}Y = 3 - M^2,$ which shows $M^2 \leq 3.$
Consequently, we get
\begin{eqnarray*}
\lefteqn{ P(Z \geq v) = P(UY \geq v-M) = E \left( P \left(U\geq \frac{v-M}Y \mid Y \right)\,{\bf 1}_{[Y > v-M]} \right) \nonumber } \\
&& \hspace{14em} = E\left(\left[1-\frac{v-M}Y \right]\, {\bf 1}_{[Y > v-M]}\right).
\end{eqnarray*}
As a similar relation holds for $P(Z \leq -u),$ we obtain the lemma.
\hfill$\Box$ \\

The following Lemma is needed in the alternative proof of Theorem \ref{unimodalCantelli}.

\begin{lemma}\label{upperboundconvex}
Let $Y$ be a real valued random variable with mean $\mu$ and variance $\sigma^2.$ For $\mu > 1$ the inequality
\begin{equation}\label{a7}
E\left(\frac 1{Y \vee 1} \right) \leq 1 - \frac{(\mu -1)^2}{\sigma^2 + \mu^2 -\mu}
\end{equation}
holds with equality if $Y$ takes on the values 1 and $\mu + \sigma^2/(\mu -1)$ with probabilities
\begin{equation}\label{a8}
P(Y=1) = \frac{\sigma^2}{\sigma^2 +(\mu - 1)^2},\quad P\left( Y = \mu + \frac {\sigma^2}{\mu - 1} \right) = \frac{(\mu - 1)^2}{\sigma^2 +(\mu - 1)^2}.
\end{equation}
For $\mu < 1$ the upper bound to $E(1/(Y \vee 1))$ is 1 with equality if $Y$ satisfies e.g. (\ref{a8}).
For $\mu = 1$ the sharp upper bound to $E(1/(Y \vee 1))$ is 1 too, but with $\sigma^2 > 0$, it cannot be attained, only approached by e.g. the Bernoulli random variable $Y_\epsilon$ defined by
\begin{eqnarray*}
P(Y_\epsilon =1 + \epsilon) = \frac{\sigma^2}{\sigma^2 + \epsilon^2},\quad P\left( Y_\epsilon = 1 - \frac {\sigma^2}{\epsilon} \right) = \frac{\epsilon^2}{\sigma^2 + \epsilon^2}, \quad \epsilon >0.
\end{eqnarray*}
Indeed $E(1/(Y_\epsilon \vee 1)) = 1- \sigma^2 \epsilon/((\sigma^2 + \epsilon^2)(1 + \epsilon)) \uparrow 1$ holds as $\epsilon \downarrow 0$.
\end{lemma}
{\bf Proof}\\
The statements in the Lemma for $\mu \leq 1$ can be checked simply by computation. We just note that the upper bound of 1 to $E(1/(Y \vee 1))$ for $\mu =1$ and $\sigma^2 >0$ cannot be attained as $Y$ has to have mass in $(1, \infty)$ then, where the function $y \mapsto 1/(y \vee 1)$ takes values strictly less than 1 only.

It remains to consider the case of prime interest, $\mu > 1.$
Observe that for any small $\varepsilon >0$ we have
\begin{eqnarray*}
E \left( \frac 1{Y \vee 1} {\bf 1}_{\left( \frac 1\varepsilon , \infty \right)}(Y) \right) < \varepsilon.
\end{eqnarray*}
Keeping this in mind we consider the conditional distribution of a random variable $Y$ with mean $\mu >1$ and variance $\sigma^2$ given $Y \leq 1/\varepsilon$.
Let $Y_\varepsilon$ be a random variable with this distribution on $(-\infty, 1/\varepsilon]$ and denote its mean and variance by $\mu_\varepsilon$ and $\sigma_\varepsilon^2$, respectively.
Note that by monotone convergence $\mu_\varepsilon \to \mu$ and $\sigma_\varepsilon^2 \to \sigma^2$ hold as $\varepsilon \downarrow 0$.
So, for sufficiently small $\varepsilon$ we have $\mu_\varepsilon >1$.
Therefore and in view of
\begin{eqnarray*}
E \left( \frac 1{Y \vee 1} \right) \leq E \left( \frac 1{Y_\varepsilon \vee 1} \right) P \left(Y \leq \frac 1 \varepsilon \right) + \varepsilon
 \leq E \left( \frac 1{Y_\varepsilon \vee 1} \right) + \varepsilon.
\end{eqnarray*}
it suffices to prove (\ref{a7}) with $Y, \mu$ and $\sigma^2$ replaced by $Y_\varepsilon, \mu_\varepsilon$ and $\sigma_\varepsilon^2$, respectively, for $\varepsilon$ sufficiently small.

If $Y_\varepsilon$ has positive mass within $(-\infty,1)$ we move this mass to 1, which increases its mean and decreases its variance, but leaves
$E(1/(Y_\varepsilon \vee 1))$ unaltered.
This increase of the mean can be compensated by moving mass of $Y_\varepsilon$ within $(\mu, 1/\varepsilon]$ towards $\mu$, thus decreasing the variance further and increasing $E(1/(Y_\varepsilon \vee 1))$.
This argument shows that we may assume that $Y_\varepsilon$ has all its mass within the interval $[1, 1/\varepsilon]$.

Let ${\cal F}(\varepsilon,\mu_\varepsilon,\tau^2)$ be the class of distribution functions with mean $\mu_\varepsilon>1$, variance $0<\tau^2 \leq \sigma_\varepsilon^2$, and all mass within $[1,1/\varepsilon]$.
This class is tight and hence by Prohorov's theorem relatively compact.
In fact it is compact.
Indeed, since any distribution $F$ in the class satisfies $F(1-)=0, F(1/\varepsilon)=1$, any limit point of a sequence within
${\cal F}(\varepsilon,\mu_\varepsilon,\tau^2)$ is a proper distribution function, and since linear and quadratic functions are continuous and bounded on $[1,1/\varepsilon]$, such a limit distribution function has mean $\mu_\varepsilon$ and variance $\tau^2$.
This compactness implies that the supremum $S(\varepsilon, \mu_\varepsilon,\tau^2)$ over ${\cal F}(\varepsilon,\mu_\varepsilon,\tau^2)$ of $E(1/(Y_\varepsilon \vee 1))$ is attained at some $F_{\varepsilon,0} \in {\cal F}(\varepsilon,\mu_\varepsilon,\tau^2).$

For $h$ a real valued, measurable, bounded function satisfying $\int h dF_{\varepsilon,0} = 0$ we may define the distribution function $F_{\varepsilon,t}$ as the one having density $1+t h \geq 0$ with respect to $F_{\varepsilon,0},$ provided $\lvert t \rvert  $ is sufficiently small.
In this way we obtain a path of distribution functions through $F_{\varepsilon,0}.$
This path is contained in ${\cal F}(\varepsilon,\mu_\varepsilon,\tau^2)$ if and only if
\begin{eqnarray*}
&& \mu_\varepsilon = \int y dF_{\varepsilon,t}(y) = \mu_\varepsilon + t \int yh(y) dF_{\varepsilon,0}(y)\ \ {\rm and} \\
&& \tau^2 + \mu_\varepsilon^2 = \int y^2 dF_{\varepsilon,t}(y) = \tau^2 + \mu_\varepsilon^2 + t \int y^2 h(y) dF_{\varepsilon,0}(y) \nonumber
\end{eqnarray*}
hold, which for $t \neq 0$ is equivalent to
\begin{eqnarray*}
\int yh(y) dF_{\varepsilon,0}(y) = \int y^2 h(y) dF_{\varepsilon,0}(y) =0.
\end{eqnarray*}
Let $\cal H$ denote the set of bounded elements $h$ in the Hilbert space $L_2(F_{\varepsilon,0})$ that are orthogonal to the quadratic polynomials. This means that the orthocomplement ${\bar {\cal H}}^\perp$ within $L_2(F_{\varepsilon,0})$ of the closed linear span $\bar {\cal H}$ of $\cal H$ consists of the set of quadratic polynomials. We have seen that each $h \in {\cal H}$ defines a smooth path through $F_{\varepsilon,0}.$
For such a path we have
\begin{eqnarray*}
\int \frac 1y dF_{\varepsilon,t}(y) = \int  \frac 1y dF_{\varepsilon,0}(y) + t \int  \frac 1y h(y) dF_{\varepsilon,0}(y).
\end{eqnarray*}
As the supremum $S(\varepsilon, \mu_\varepsilon, \tau^2)$ is attained at $F_{\varepsilon,0},$ this implies for all $h \in {\cal H}$
\begin{eqnarray*}
\int  \frac 1y h(y) dF_{\varepsilon,0}(y) = 0
\end{eqnarray*}
and we may conclude that $1/y$ belongs to ${\bar {\cal H}}^\perp.$
Consequently, $1/y$ is a quadratic polynomial $F_{\varepsilon,0}$ almost surely.

This means that there exist $p_0,p_1$ and $p_2$ such that
\begin{eqnarray*}
p_0 + p_1 y + p_2 y^2 = 1/y {\rm ~or~} p_0 y + p_1 y^2 + p_2 y^3 = 1
\end{eqnarray*}
holds for $F_{\varepsilon,0}$ almost all $y \in [1,1/\varepsilon]$.
As cubic polynomials have at most 3 roots, we shall determine the supremum of $E(1/Y_\varepsilon)$ for random variables
\begin{equation}\label{a17}
Y_\varepsilon = \left\{
\begin{array}{rcl}
a & & p \\
b & \mbox{with probability} & q \\
c & & 1-p-q \\
\end{array}\right.
\end{equation}
with $a < b < c \in [1,1/\varepsilon], E(Y_\varepsilon) = \mu_\varepsilon >1$ and ${\rm var}(Y_\varepsilon) = \tau^2 \leq \sigma_\varepsilon^2$.
We consider the following cases.
\begin{itemize}
\item[1 atom]
In this case $Y_\varepsilon$ has to be degenerate at $\mu_\varepsilon$ and $E(1/Y_\varepsilon) = 1/\mu_\varepsilon$ holds, which is a value smaller than the upper bound from (\ref{a7}) with $\mu$ and $\sigma^2$ replaced by $\mu_\varepsilon$ and $\tau^2$ or $\sigma_\varepsilon^2$.
\item[2 atoms]
Here we may take $q=1-p$.
Solving
\begin{eqnarray*}
\mu_\varepsilon = ap + b(1-p),\quad \tau^2 + \mu_\varepsilon^2 = a^2 p + b^2(1-p)
\end{eqnarray*}
for $b$ and $p$ we obtain $b = \mu_\varepsilon + \tau^2 / (\mu_\varepsilon -a)$ and $p = \tau^2 /( \tau^2 + (\mu_\varepsilon -a)^2),$
and hence
\begin{eqnarray*}
E \left( \frac 1{Y_\varepsilon} \right) = \frac 1a p + \frac 1b (1-p)
= \frac 1{\mu_\varepsilon} \left[1 + \frac {\tau^2}{a \left\{ \mu_\varepsilon + \tau^2/(\mu_\varepsilon -a) \right\} } \right]
\end{eqnarray*}
with $1 \leq a \leq \mu_\varepsilon$.
As the denominator in this expression is increasing in $a$, the supremum is attained at $a=1.$
The resulting distribution of $Y_\varepsilon$ has two atoms, one of which is at 1, and satisfies (\ref{a8}) and the right hand side of (\ref{a7}) with $\mu$ and $\sigma^2$ replaced by $\mu_\varepsilon$ and $\tau^2.$
In view of $\tau^2 \leq \sigma_\varepsilon^2$ the $\varepsilon$ version of the upper bound in (\ref{a7}) is valid for the two atoms case we are considering here, and is attained by the $\varepsilon$ version of the distribution from (\ref{a8}).
\item[3 atoms]
In this most general case of (\ref{a17}) consider the conditional expectation $E(1/Y_\varepsilon \mid Y_\varepsilon \neq c)$ given the values of $E(Y_\varepsilon \mid Y_\varepsilon \neq c)$ and ${\rm var}(Y_\varepsilon \mid Y_\varepsilon \neq c).$ This brings us back to the preceding case with two atoms, and we may conclude that the supremum is attained for $a=1.$
Solving
\begin{eqnarray*}
\mu_\varepsilon = p + bq + c(1-p-q), \quad \tau^2 + \mu_\varepsilon^2 = p + b^2 q + c^2 (1-p-q)
\end{eqnarray*}
for $p$ and $q$ we obtain
\begin{eqnarray}\label{a20}
&& p = \frac{\tau^2 +(b-\mu_\varepsilon)(c-\mu_\varepsilon)}{(b-1)(c-1)},\quad
q = \frac{(\mu_\varepsilon -1)(c-\mu_\varepsilon) - \tau^2}{(b-1)(c-b)} \nonumber \\
&& E \left( \frac 1Y_\varepsilon \right) = 1 - \frac {\mu_\varepsilon -1}c + \left[ \tau^2 -(c- \mu_\varepsilon)(\mu_\varepsilon -1) \right] \frac 1{bc}.
\end{eqnarray}
Since the variance of a random variable on a bounded interval is largest if all its mass is at the boundary of the interval as in the Bernoulli case, we have $\tau^2 \leq (c- \mu_\varepsilon)(\mu_\varepsilon -1)$.
Consequently, the maximum of the last expression in (\ref{a20}) over $b \in [1,c]$ is attained at $b=c$, which brings us back to the two atoms case.
\end{itemize}
As argued above, taking limits with $\varepsilon \downarrow 0$ completes the proof of (\ref{a7}).
Simple computation shows that the distribution from (\ref{a8}) attains this bound.
\hfill$\Box$ \\

In the proof of Theorem \ref{generalGauss} we will need to solve the following cubic equation.

\begin{lemma}\label{cubicequation}
With $r > 0$ the cubic equation
\begin{equation}\label{cubicequation1}
z^3 + 3 z^2 - 3r z -r =0
\end{equation}
has exactly one positive root, namely
\begin{eqnarray*}
z_r = 2 \sqrt{1+r} \cos \left( \frac 13 \left[\pi -\arctan \left(\sqrt{r}\right)\right] \right) - 1
\end{eqnarray*}
with, in particular, $z_1 =1.$
\end{lemma}
{\bf Proof}\\
By Vieta's method the equation (\ref{cubicequation1}) becomes
\begin{eqnarray*}
\left(w^3 \right)^2 +  2(1+r)w^3 + (1+r)^3 =0
\end{eqnarray*}
via the  substitution $z=w +(1+r)w^{-1} -1.$ This leads to
\begin{eqnarray*}
w^3 = -(1+r)(1 \pm \sqrt{r}\,i) = -(1+r)^{3/2} \exp\left(2k\pi i \pm \arctan \left(\sqrt{r}\right)\right)
\end{eqnarray*}
and hence to
\begin{eqnarray*}
\lefteqn{z = - 1 - \sqrt{1+r} \nonumber } \\
&& \left[ \exp\left(\frac 13 \left[2k\pi \pm \arctan \left(\sqrt{r}\right)\right] i \right) +
 \exp\left(-\frac 13 \left[2k\pi \pm \arctan \left(\sqrt{r}\right)\right] i \right)\right] \nonumber \\
&& = - 2\sqrt{1+r}\, \cos \left(\frac 13 \left[2k\pi + \arctan \left(\sqrt{r}\right)\right]\right) -1 \\
&& = 2\sqrt{1+r}\, \cos \left(\pi - \frac 13 \left[2k\pi + \arctan \left(\sqrt{r}\right)\right]\right) - 1 = 2\sqrt{1+r}\, \cos (\alpha) -1 \nonumber
\end{eqnarray*}
with $k$ an integer and $\alpha$ equal to
\begin{eqnarray*}
&& \alpha_0 = \frac 13 \left[ 3\pi - \arctan \left( \sqrt{r} \right)\right] \in (5\pi/6, \pi), \nonumber \\
&& \alpha_1 = \frac 13 \left[ \pi - \arctan \left( \sqrt{r} \right)\right] \in (\pi /6, \pi/3), {\rm ~ or} \\
&& \alpha_2 = \frac 13 \left[ - \pi - \arctan \left( \sqrt{r} \right)\right] \in (-\pi /2, -\pi/3). \nonumber
\end{eqnarray*}
In view of
\begin{eqnarray*}
\cos(\alpha_0) < -\sqrt{3}/2,\ \cos(\alpha_1) > 1/2,
\end{eqnarray*}
the claim has been proved once
\begin{equation}\label{cubicequation5}
2 \sqrt{1+r}\, \cos(-\alpha_2) - 1 < 0
\end{equation}
has been shown. Writing $\psi = \pi/6 -\arctan(\sqrt{r})/3$ with $\psi \in (0, \pi/6)$ we have $r = \cot^2(3\psi)$ and we see that (\ref{cubicequation5}) holds if and only if
\begin{eqnarray*}
\frac 2{\sin(3\psi)} \cos(\pi/2 - \psi) -1 < 0
\end{eqnarray*}
or equivalently
\begin{eqnarray*}
2 \sin(\psi) - \sin(3\psi) = 2\sin(\psi) - 3\sin(\psi) + 4 \sin^3(\psi) = \sin(\psi)[4\sin^2(\psi) -1] < 0
\end{eqnarray*}
holds. Because of $\sin(\psi) \in (0, 1/2)$ this is the case.
\hfill$\Box$ \\

\subsection{Alternative Proof of Theorem \ref{unimodalCantelli}}\label{proofunimodalCantelli}

Here we present an alternative to the proof of Theorem \ref{unimodalCantelli} as given in Section \ref{unimodal}. \\
\noindent
{\bf Preparatory Part of the Proof}\\
Let ${\cal Z}_M$ and ${\cal Y}_M$ be the classes of distributions as defined in Lemma \ref{unimodalCantelliPsi}.
First we consider the case that the location of the mode satisfies $M < v.$
By Lemma \ref{unimodalCantelliPsi} with $u=\infty$ we obtain
\begin{equation}\label{Ppsi}
\sup_{Z \in {\cal Z}_M} P(Z\geq v) = \sup_{Y \in {\cal Y}_M} E\psi_M(Y)
\end{equation}
with $\psi_M$ defined by
\begin{eqnarray*}
\psi_M(y) = \left[ 1-\frac{v-M}y \right] {\bf 1}_{[y>v-M]}.
\end{eqnarray*}
With $Y \in {\cal Y}_M$ we define the Bernoulli random variable $Y_2$ by
\begin{eqnarray*}
Y_2= \mu_-\,{\bf 1}_{[Y<v-M]}\ +\ \mu_+\,{\bf 1}_{[Y\geq v-M]},
\end{eqnarray*}
with $\mu_- = E(Y \mid Y<v-M)$ and $\mu_+ = E(Y \mid Y\geq v-M)$.
Note $EY_2=EY$ and $EY_2^2 \leq EY^2.$ As $E(Y \mid Y<v-M) <
v-M\leq E(Y \mid Y\geq v-M)$ holds and $\psi_M$ is concave
on $[v-M,\infty)$ and vanishes elsewhere, we have by Jensen's
inequality
\begin{eqnarray*}
E\psi_M(Y) & = & E(\psi_M(Y) \mid Y\geq v-M)P(Y\geq v-M) \\
&\leq & \psi_M(E(Y \mid Y\geq v-M))P(Y\geq v-M) = E\psi_M(Y_{2}).
\end{eqnarray*}
By adding a positive amount to $\mu_+$ and subtracting from $\mu_-$, if necessary, we can force the
Bernoulli random variable $Y_2$ to have variance $3-M^2$ while we maintain its mean at $-2M$ and possibly increase $E\psi_M(Y_{2}),$ as $\psi_M$ is increasing on $[v-M, \infty).$
We have shown that the supremum at the right hand side of (\ref{Ppsi}) is attained by a Bernoulli random variable
\begin{eqnarray*}
Y_{2} = \left\{
\begin{array}{rcl}
a & & 1 - p\\
   &\mbox{with probability}& \\
b & & p
\end{array}\right.
\end{eqnarray*}
with $a< v - M \leq b,$ and $0\leq p\leq 1$.  In view of $EY_{2}=-2M$
and $EY_{2}^2=3(1+M^2)$ we obtain
\begin{equation}\label{panda}
p=\frac{3 - M^2}{(b+2M)^2 + 3 - M^2},\quad a=-\,\frac{2bM+3(1+M^2)}{b+2M}.
\end{equation}
Writing
\begin{equation}\label{defb}
b=(v-M)(1+x),\quad x \geq 0,
\end{equation}
we see that the suprema from (\ref{Ppsi}) equal
\begin{eqnarray}\label{supg}
\lefteqn{\sup_{b \geq v-M >0} \left[1-\frac{v-M}{b}\right]p = \sup_{x \geq 0,M < v}\
\frac x{1+x}\ \frac{3 - M^2}{\{ (1+x)v + (1-x)M \}^2 + 3 - M^2} \nonumber}\\
&& \nonumber \\
&& \hspace{11em} = \sup_{x \geq 0} \frac x{1+x} \left[1 + \inf_{M<v} \frac{\{ (1+x)v + (1-x)M \}^2}{3 - M^2} \right]^{-1} \nonumber \\
&& \nonumber \\
&& \hspace{11em} = \sup_{x \geq 0} \frac x{1+x} \left[1 + \inf_{M<v} \chi(M \mid x) \right]^{-1},
\end{eqnarray}
where it is understood that $M$ is restricted by $M^2 < 3.$ Given $x \geq 0$ we study $\chi(M \mid x)$ as a function of $M$ and note that this function has vertical asymptotes at $M= - {\sqrt 3}$ and $M={\sqrt 3}$ with
\begin{equation}\label{asymptotes}
\lim_{M \downarrow -{\sqrt 3}} \chi(M \mid x) = \lim_{M \uparrow {\sqrt 3}} \chi(M \mid x) = \infty.
\end{equation}
It has two stationary points, namely
\begin{equation}\label{stationarypoints}
M_1(x)= -\ \frac{(1+x)v}{1-x}\ {\rm and}\ M_2(x) = -\ \frac{3(1-x)}{(1+x)v}.
\end{equation}
In view of (\ref{asymptotes}) but also because of $M_1(x) M_2(x) =3$ exactly one of these stationary points belongs to the interval $(-{\sqrt 3}, {\sqrt 3}).$
Indeed, $M_1^2(x) < 3$ holds with $\chi(M_1(x) \mid x)=0,$ if and only if $v < {\sqrt 3}$ holds and
\begin{equation}\label{xforM1}
x< \frac{{\sqrt 3} -v}{{\sqrt 3} + v}\ {\rm or}\  \frac{{\sqrt 3} + v}{{\sqrt 3} -v} < x
\end{equation}
is valid. On the other hand,
$M_2^2(x) < 3$ holds if and only if $v \geq {\sqrt 3}$ holds or $v < {\sqrt 3}$ is valid with
\begin{equation}\label{xforM2}
\frac{{\sqrt 3} -v}{{\sqrt 3} + v} < x < \frac{{\sqrt 3} + v}{{\sqrt 3} -v}.
\end{equation}
Furthermore, we have
\begin{equation}\label{valuechiM2}
\chi(M_2(x) \mid x) = \frac 13 v^2(1+x)^2 - (1-x)^2,
\end{equation}
which is positive for $M_2^2(x) < 3.$
For $M < v$ we have collected some properties of the supremum in (\ref{Ppsi}) that we need in the proofs below.
\hfill$\Box$\\

\noindent
{\bf Proof of (\ref{unimodalCantellia}) for $v > {\sqrt 3}$}\\
As the location of the mode cannot be outside the interval $[-{\sqrt 3}, {\sqrt 3}],$
we are in the situation of the {\bf preparatory part of the proof} with $M<v.$ Consequently, (\ref{Ppsi}), (\ref{supg}), and (\ref{valuechiM2}) yield
\begin{equation}\label{supchi}
P(Z \geq v) \leq \sup_{x \geq 0} \frac x{(1+x) \left[1 + \frac 13 v^2(1+x)^2 - (1-x)^2 \right]}.
\end{equation}
As the denominator is strictly positive for nonnegative $x$ and the function itself vanishes at 0 and infinity, the supremum is attained at one of its stationary points on the positive half line. Straightforward computations show that they satisfy the equation
\begin{eqnarray*}
2(v^2-3)x^3 + 3(v^2 +1)x^2 - v^2 = 0,
\end{eqnarray*}
which can be rewritten as
\begin{eqnarray*}
(2x-1)\left[(v^2 -3)x^2 + 2v^2 x + v^2 \right]=0.
\end{eqnarray*}
As the quadratic factor has negative roots $-v/(v + {\sqrt 3})$ and $-v/(v-{\sqrt 3}),$ the only positive stationary point is at $x=1/2,$ and (\ref{supchi}) with $x=1/2$ yields the second bound from (\ref{unimodalCantellia}), namely $4/(9(1+v^2)).$

At $x=1/2$ equations (\ref{stationarypoints}), (\ref{defb}), and (\ref{panda}) yield $M_2(1/2)= -1/v, b = 3(1+v^2)/(2v), p = 4/(3 + 3v^2),$ and $a=0.$ These values determine the distribution of the random variable $Z$ that attains equality in (\ref{unimodalCantellia}). An alternative description of this distribution is given in the statement of the Theorem.
\hfill$\Box$\\

\noindent
{\bf Proof of (\ref{unimodalCantellia}) for $0 \leq v \leq {\sqrt 3}$}\\
In view of $M^2 \leq 3$ we discern three cases, namely $-{\sqrt 3} \leq M < v, M=v,$ and $v < M \leq {\sqrt 3}.$
First, we study the case $-{\sqrt 3} \leq M < v,$  for which the {\bf preparatory part of the proof} can be applied. We partition the nonnegative half line into the intervals
\begin{eqnarray*}
&& I_1 = \left[0,\frac{{\sqrt 3} - v}{{\sqrt 3} +v}\right),\quad I_2 =\left[\frac{{\sqrt 3} - v}{{\sqrt 3} +v}, \frac{3 + v^2}{3-v^2}\right], \nonumber \\
&& I_3 = \left(\frac{3 + v^2}{3-v^2}, \frac{{\sqrt 3} + v}{{\sqrt 3} - v}\right),\quad I_4 = \left[\frac{{\sqrt 3} + v}{{\sqrt 3} - v}, \infty \right).
\end{eqnarray*}
Straightforward computation, (\ref{stationarypoints}), (\ref{xforM1}), and (\ref{xforM2}) show
\begin{eqnarray*}
&& x \in I_1 \Rightarrow -{\sqrt 3} \leq M_1(x) \leq v, \quad x \in I_2 \Rightarrow -{\sqrt 3} \leq M_2(x) \leq v, \nonumber \\
&& \\
&& x \in I_3 \Rightarrow v < M_2(x) \leq {\sqrt 3}, \quad x \in I_4 \Rightarrow  v < M_1(x) \leq {\sqrt 3}, \nonumber
\end{eqnarray*}
which by (\ref{Ppsi}) and (\ref{supg}) yields
\begin{eqnarray}\label{maxs}
\lefteqn{\sup_{M<v} P(Z \geq v) \leq \max \left\{\sup_{x \in I_1} \frac x{1+x},\ \sup_{x \in I_2} \frac x{1+x}\ \frac 1{1+\chi(M_2(x) \mid x)},\right. \nonumber } \\
&& \qquad \qquad \qquad \qquad \left. \sup_{x \in I_3} \frac x{1+x}\ \frac 1{1+\chi(v \mid x)},\ \sup_{x \in I_4} \frac x{1+x}\ \frac 1{1+\chi(v \mid  x)} \right\} \nonumber \\
&& = \max \left\{\frac{{\sqrt 3} - v}{2{\sqrt 3}},\ \sup_{x \in I_2} \frac x{1+x}\ \frac 1{1+\chi(M_2(x) \mid x)},\right. \\
&& \qquad \qquad \left. \sup_{x \in I_3} \frac x{1+x}\ \frac {3-v^2}{3(1+v^2)},\ \sup_{x \in I_4} \frac x{1+x}\ \frac {3-v^2}{3(1+v^2)} \right\} \nonumber \\
&& = \max \left\{ \frac {3-v^2}{3(1+v^2)},\ \sup_{x \in I_2} \frac x{(1+x) \left[1 + \frac 13 v^2(1+x)^2 - (1-x)^2 \right]} \right\}. \nonumber
\end{eqnarray}
As shown in the {\bf Proof of (\ref{unimodalCantellia}) for $v > {\sqrt 3}$} the supremum at the right hand side of (\ref{maxs}) equals $4/(9(1+v^2)),$ provided $1/2 \in I_2$ holds. Now, 1/2 belongs to $I_2$ if and only if $v \geq \sqrt{1/3}$ holds. Furthermore, $4/(9(1+v^2))$ equals at least $(3-v^2)/(3(1+v^2))$ if and only if $v \geq \sqrt{5/3}$ holds. We may conclude
\begin{equation}\label{maxMsv}
\sup_{M<v} P(Z \geq v) \leq \max \left\{\frac {3-v^2}{3(1+v^2)},\ \frac 4{9(1+v^2)} \right\}.
\end{equation}

For standardized random variables $Z=M + UY$ with $v < M \leq {\sqrt 3}$ we have
\begin{eqnarray*}
Z \geq v \Leftrightarrow \frac {UY}{v-M} \leq 1 \Leftrightarrow U \left( \left( \frac{-Y}{M-v} \right) \vee 1 \right) \leq 1.
\end{eqnarray*}
Consequently, Lemma \ref{upperboundconvex} with $\mu= 2M/(M-v)$ and $\sigma^2 = (3-M^2)/(M-v)^2$ yields
\begin{eqnarray}\label{maxvsM}
\lefteqn{ \sup_{M>v} P(Z \geq v) = \sup_{M>v} E\left( \frac 1{(-Y/(M-v)) \vee 1} \right) \nonumber }\\
&& \leq \sup_{M>v} \frac {3-v^2}{3 + 2vM + M^2} = \frac {3-v^2}{3(1+v^2)}.
\end{eqnarray}
The case $M=v$ may be settled by noting that (\ref{maxvsM}) implies
\begin{eqnarray*}
P(Z \geq v) = \lim_{\epsilon \downarrow 0} P(Z \geq v - \epsilon) \leq \frac {3-v^2}{3(1+v^2)}.
\end{eqnarray*}
Combining this with (\ref{maxMsv}) and (\ref{maxvsM}) we complete the proof of (\ref{unimodalCantellia}). By straightforward computation one may verify that the distributions mentioned in the statement of the Theorem attain the bounds in (\ref{unimodalCantellia}).
\hfill$\Box$\\

\subsection{Proof of Theorem \ref{generalGauss}}\label{proofgeneralGauss}
Finally, we present our proof of Theorem \ref{generalGauss} as given in Section \ref{unimodal}. \\
By Khintchine's Lemma \ref{Khintchinelemma} the standardized unimodal random variable $Z$ may be represented as $Z=M + UY$ with $U$ uniformly distributed on the unit interval and independent of the random variable $Y.$
As $Z$ is standardized, $Y$ has to satisfy $EY=-2M$ and $EY^2=3(1+M^2)$ with $M$ the location of the mode of $Z.$
It follows that the variance of $Y$ equals $3-M^2,$ and hence that $\lvert M \rvert   \leq \sqrt 3$ holds.
As $u$ and $v$ are both at least as large as $\sqrt 3$, we have $-u \leq M \leq v$.
Hence, Khintchine's representation and Lemma \ref{unimodalCantelliPsi} yield
\begin{equation}\label{s1a}
P(Z\leq -u\ \  \mbox{\rm or}\ \  Z \geq v) = E \Psi_M(Y),
\end{equation}
where the function $\Psi_M$ is given by
\begin{eqnarray*}
\Psi_M(y)=\left[1+\frac{u+M}y \right] {\bf 1}_{[y<-u-M]}+ \left[1-\frac{v-M}y \right] {\bf 1}_{[v-M < y]}.
\end{eqnarray*}

We define the random variable $Y_1$ by
\begin{eqnarray*}
Y_1= \mu_-\,{\bf 1}_{[Y<-u-M]} + \mu_0\,{\bf 1}_{[-u-M \leq Y \leq v-M]} + \mu_+\,{\bf 1}_{[v-M < Y]},
\end{eqnarray*}
with $\mu_- = E(Y \mid Y< -u-M),\ \mu_0 = E(Y \mid -u-M \leq Y \leq v-M),$ and $\mu_+ = E(Y \mid Y > v-M)$.
Note $EY_1=EY$ and $EY_1^2 \leq EY^2.$

As $\mu_- < -u-M \leq \mu_0 \leq v-M < \mu_+$ holds and $\psi_M$ is concave
on $(-\infty, -u-M)$ and on $(v-M,\infty)$ and vanishes elsewhere, we have by Jensen's
inequality
\begin{eqnarray*}
\lefteqn{E\Psi_M(Y) \nonumber }\\
&& = [E(\Psi_M(Y) \mid  Y < -u-M)]\,P(Y < -u-M) \nonumber \\
&& \qquad \qquad + [E(\Psi_M(Y) \mid v-M < Y)]\, P(v-M < Y) \nonumber \\
&& \leq \Psi_M(E(Y \mid Y < -u-M))P(Y < -u-M) \\
&& \qquad \qquad + \Psi_M(E(Y \mid v-M < Y))P(v-M < Y) \nonumber \\
&& = \Psi_M(\mu_-)P(Y_1 = \mu_-) + \Psi_M(\mu_+)P(Y_1 = \mu_+) = E\Psi_M(Y_1). \nonumber
\end{eqnarray*}
If necessary, by subtracting a positive value from $\mu_-$ (or $\mu_0$ if the mass at $\mu_-$ vanishes) and adding to $\mu_+$ (or $\mu_0$ if the mass at $\mu_+$ vanishes), the random variable $Y_1$ can be forced to satisfy $E(Y_1^2)=3(1+M^2)$ while increasing $E\Psi_M(Y_1)$ and maintaining $EY_1=-2M.$
This mechanism does not work if $Y_1$ is degenerate at $-2M$.
However, in view of $-u-M \leq -2M \leq v-M$ we then have $E\Psi_M(Y_1)=0$, which is not an upper bound to (\ref{s1a}) whatever the values of $u$ and $v$ are.

We have shown that the supremum of the probability in (\ref{boundgeneralGauss}) is attained by a random variable $Z=M + UY$ as above, where $Y$ is discrete with three atoms, namely
\begin{eqnarray*}
Y = \left\{
\begin{array}{rcl}
-a & & p \\
b  &\mbox{with probability}& 1-p-q \\
c & & q
\end{array}\right.
\end{eqnarray*}
with $-a \leq -u-M \leq b \leq v-M \leq c, \lvert M \rvert   \leq \sqrt{3},$ and
\begin{equation}\label{G0a}
P(Z\leq -u\ \  \mbox{\rm or}\ \  Z \geq v) = \left(1-\frac{u+M}a \right)p + \left(1-\frac{v-M}c \right) q.
\end{equation}
The restrictions $EZ=0$ and $EZ^2=1$ imply
\begin{eqnarray}\label{G1a}
EY & = & -2M = -ap + b(1-p-q) + cq \\
EY^2 & = & 3(1+M^2) = a^2p + b^2(1-p-q) + c^2q, \nonumber
\end{eqnarray}
and hence
\begin{eqnarray}\label{pandq}
\lefteqn{p=\frac{3(1+M^2)+2M(b+c)+bc}{(a+b)(a+c)},\quad q = \frac{3(1+M^2)+2M(b-a)-ab}{(c-b)(a+c)}, \nonumber } \\
&& \hspace{3em} p+q=\frac{3(1+M^2)+2M(c-a)+b(c-a-b)}{(a+b)(c-b)}. \hspace{3em}
\end{eqnarray}
Writing
\begin{equation}\label{zetaeta}
\zeta = 3(1+M^2) + 2Mb, \quad \eta = 2M+b,
\end{equation}
we note
\begin{eqnarray}\label{pqpositive}
p \geq 0 & \Leftrightarrow & bc+2M(b+c)+3(1+M^2) = \zeta + \eta c \geq 0, \nonumber \\
q \geq 0 & \Leftrightarrow & ab+2M(a-b)-3(1+M^2) = -\zeta + \eta a \leq 0, \\
p+q \leq 1 & \Leftrightarrow & ac+2M(a-c) -3(1+M^2) \geq 0, \nonumber
\end{eqnarray}
and define the set
\begin{eqnarray*}
\lefteqn{ {\cal A} = \left\{(a,b,c,M) \, \mid \, -a \leq -u-M \leq b \leq v-M \leq c, \lvert M \rvert   \leq \sqrt{3},\right. \nonumber } \\
&& \left. \hspace{11em} p\geq 0, q \geq 0, p+q \leq 1 \right\}.
\end{eqnarray*}
So far, we have seen that
\begin{eqnarray}\label{maximumoverA}
\lefteqn{P(Z\leq -u\ \  \mbox{\rm or}\ \  Z \geq v) } \\
&& \leq \sup_{{\cal A}} \left[ \left( 1 - \frac{u+M}a \right)\frac{\zeta + \eta c}{(a+b)(a+c)}
+ \left(1 - \frac{v-M}c \right)\frac{\zeta - \eta a}{(c-b)(a+c)} \right] \nonumber
\end{eqnarray}
holds.
Note that $u+M$ and $v-M$ are positive on $\cal A$ in view of $-u \leq M \leq v$.

We shall prove that for $\sqrt{3} \leq v \leq u \leq v+2/v$ this supremum is attained at a stationary point within $\cal A$ of the function in (\ref{maximumoverA}) and that for $v + 2/v \leq u$ it is attained at a point on the boundary of $\cal A$.
To this end we shall show first that at the boundary of $\cal A$ the function in (\ref{maximumoverA}) cannot attain a value larger than the second bound given in (\ref{boundgeneralGauss}).
With $\bar{\cal A}$ denoting the closure of $\cal A$ we see that the boundary $\partial A$ of $\cal A$ is the union of the following sets
\begin{eqnarray*}
&& {\cal A}_1 = \{(a,b,c,M) \in {\cal A} \, \mid \, a=u+M \}, \nonumber \\
&& {\cal A}_2 = \{(a,b,c,M) \in {\cal A} \, \mid \, c=v-M \}, \nonumber \\
&& {\cal A}_3 = \{(a,b,c,M) \in {\cal A} \, \mid \, \zeta - \eta a = 0 \} = \{(a,b,c,M) \in {\cal A} \, \mid \, q=0 \}, \nonumber \\
&& {\cal A}_4 = \{(a,b,c,M) \in {\cal A} \, \mid \, \zeta + \eta c = 0 \} = \{(a,b,c,M) \in {\cal A} \, \mid \, p=0 \}, \nonumber \\
&& {\cal A}_5 = \{(a,b,c,M) \in {\cal A} \, \mid \, p+q=1 \}, \\
&& {\cal A}_6 = \{(a,b,c,M) \in \bar{\cal A} \, \mid \, \lvert M \rvert   = \sqrt{3} \}, \nonumber \\
&& {\cal A}_7 = \{(a,b,c,M) \in \bar{\cal A} \, \mid \, b=-a \}, \nonumber \\
&& {\cal A}_8 = \{(a,b,c,M) \in \bar{\cal A} \, \mid \, b=c \}. \nonumber
\end{eqnarray*}
We treat these boundary subsets as follows.
\begin{itemize}
\item[${\cal A}_1$]
With $a=u+M$ we have
\begin{equation}\label{B1}
P(Z\leq -u\ \  \mbox{\rm or}\ \  Z \geq v) = \left(1 - \frac{v-M}c \right)q = P(Z \geq v),
\end{equation}
which by Theorem \ref{unimodalCantelli} is bounded by $4/(9(1+v^2))$ provided $v^2 \geq 5/3$ holds.
By differentiation one observes that the function $$u \mapsto \frac{4+(u-v)^2}{(u+v)^2}$$ is decreasing if and only if $u \leq v+2/v$ holds, and hence it has minimum value $1/(1+v^2)$.
So, the first bound from (\ref{boundgeneralGauss}) equals at least the bound $4/(9(1+v^2))$ from Theorem \ref{unimodalCantelli}; cf. the alternative proof in the Appendix.

\item[${\cal A}_2$]
By symmetry an analogous argument holds for ${\cal A}_2$ as for ${\cal A}_1$.

\item[${\cal A}_3$]
By symmetry an analogous argument holds for ${\cal A}_3$ as for ${\cal A}_4$.

\item[${\cal A}_4$]
With $p=0$ we have (\ref{B1}) and the same argument as for ${\cal A}_1$ holds here.
Furthermore, the random variable $Z$ that attains the second bound from (\ref{boundgeneralGauss}), corresponds to $M=-1/v$ and
\begin{eqnarray*}
Y = \left\{
\begin{array}{rcl}
0 & & \frac {3v^2 - 1}{3(1+v^2)} \\
 &\mbox{with probability}&  \\
\frac 32 \left( v + \frac 1v \right) & & \frac 4{3(1+v^2)},
\end{array}\right.
\end{eqnarray*}
which shows that this second bound is attained within ${\cal A}_4$.

\item[${\cal A}_5$]
In view of $EZ=0$ and $EZ^2=1$, the definition of $Y$ from (\ref{G1a}) implies $EY=-2M$ and $EY^2=3(1+M^2)$, and hence $E(Y+M)^2 =3$.
With $p+q=1$ this yields
\begin{eqnarray*}
(a-M)^2 p + (M+c)^2 (1-p) =3,
\end{eqnarray*}
which means that $a-M$ and $M+c$ cannot simultaneously be larger than $\sqrt{3}$.
As both $u$ and $v$ equal at least $\sqrt{3}$, this shows that either $1-(u+M)/a \leq 0$ or $1-(v-M)/c \leq 0$ holds.
Together with (\ref{B1}) we conclude that ${\cal A}_5 \subset {\cal A}_1 \cup {\cal A}_2$ holds and that the second bound from (\ref{boundgeneralGauss}) cannot be exceeded on ${\cal A}_5$.

\item[${\cal A}_6$]
In case of $\lvert M \rvert   = \sqrt{3}$ the variance of $Y$ from (\ref{G1a}) vanishes, i.e., $Y$ is degenerate, and hence at least two of the restrictions $p=0, q=0$, and $p+q=1$ hold. Consequently, we have ${\cal A}_6 \subset {\cal A}_3 \cup {\cal A}_4$ and we see that the second bound from (\ref{boundgeneralGauss}) cannot be exceeded on ${\cal A}_6$.

\item[${\cal A}_7$]
If $b$ equals $-a$, the random variable $Y$ from (\ref{G1a}) may be viewed as a Bernoulli random variable with $p+q=1$, which implies ${\cal A}_7 \subset {\cal A}_5$.

\item[${\cal A}_8$]
If $b$ equals $c$, the random variable $Y$ from (\ref{G1a}) may be viewed as a Bernoulli random variable with $p+q=1$, which implies ${\cal A}_8 \subset {\cal A}_5$.

\end{itemize}

We conclude that for $\sqrt{3} \leq v \leq u$
\begin{eqnarray}\label{maximumoverboundaryA}
\lefteqn{ \sup_{\partial{\cal A}} \left[ \left( 1 - \frac{u+M}a \right)\frac{\zeta + \eta c}{(a+b)(a+c)}
+ \left(1 - \frac{v-M}c \right)\frac{\zeta - \eta a}{(c-b)(a+c)} \right] \nonumber } \\
&& \hspace{19em} = \frac 49\ \frac 1{1+v^2}
\end{eqnarray}
holds.
As we have shown that at the boundary of $\cal A$ the function from (\ref{maximumoverA}) cannot attain a value larger than the first bound given in (\ref{boundgeneralGauss}), we focus on the interior of $\cal A$ and determine the stationary points of the function in (\ref{maximumoverA}).
For the time being we fix $b$ and hence $\zeta$ and $\eta$ and note that the function to be maximized over $a$ and $c$ may be written as
\begin{eqnarray}\label{functiontomaximizea}
\lefteqn{ \left\{\zeta [ac(a+c) - (u+M)c(c-b) - (v-M)a(a+b)] \right. } \\
&& \left. + \eta [ac(a+c)(c-a-b) - (u+M)c^2(c-b) + (v-M)a^2(a+b)] \right\} \nonumber \\
&& \times \{a(a+b)c(c-b)(a+c)\}^{-1}. \nonumber
\end{eqnarray}
"Some" computation shows that the stationary points of this function of $a$ and $c$ are solutions of the two equations
\begin{eqnarray*}
&& [\zeta + \eta c]\ [-a^2c(a+c)^2 + (u+M)c(c-b)(3a^2 + 2a(b+c) + bc) \nonumber \\
&& \qquad \qquad \qquad \qquad \qquad + (v-M)a^2(a+b)^2] =0 \\
&& [\zeta - \eta a]\ [-ac^2(a+c)^2 + (u+M)c^2(c-b)^2 \nonumber \\
&& \qquad \qquad \qquad \qquad \qquad + (v-M)a(a+b)(3c^2 + 2c(a-b) - ab)] =0 . \nonumber
\end{eqnarray*}
Ignoring the first factors, which correspond to the boundary conditions $p=0$ and $q=0$ treated under ${\cal A}_3$ and ${\cal A}_4$ above, we obtain
\begin{eqnarray}\label{LA61a}
\lefteqn{a^2c^2(a+c)^2 }\\
&& = (u+M)c^2(c-b)(3a^2 + 2a(b+c) + bc) + (v-M)a^2(a+b)^2 c \nonumber \\
&& =  (u+M)ac^2(c-b)^2 + (v-M)a^2(a+b)(3c^2 + 2c(a-b) - ab). \nonumber
\end{eqnarray}
Dividing the second equality by $v-M$ we obtain
\begin{eqnarray*}
\frac{u+M}{v-M}\, c^2(c-b)(a+b)(3a+c) + a^2(a+b)(c-b)(-3c-a) = 0.
\end{eqnarray*}
Dividing this by $c^3(c-b)(a+b)$ and writing $a = \gamma c$ we arrive at
\begin{equation}\label{LA63a}
\gamma^3 + 3 \gamma^2 -3 \frac{u+M}{v-M} \gamma - \frac{u+M}{v-M} = 0,
\end{equation}
which by Lemma \ref{cubicequation} has exactly one positive root, namely
\begin{equation}\label{gammaa}
\gamma_M = 2 \sqrt{1+ \frac{u+M}{v-M}} \cos \left( \frac 13 \left[\pi -\arctan \left(\sqrt{\frac{u+M}{v-M}}\right)\right] \right) - 1 > 0.
\end{equation}
By a slight abuse of notation we shall denote this unique positive root by $\gamma$ too.
We conclude that $a$ and $c$ satisfy (\ref{LA61a}) and
\begin{eqnarray*}
a = \gamma c
\end{eqnarray*}
holds with $a$ and $c$ depending on $b$ and $M$, and with $\gamma$ depending on $M$ only.

Note that $\zeta$ and $\eta$ defined in (\ref{zetaeta}) depend on $b$.
Straightforward computation shows that the derivative of (\ref{functiontomaximizea}) with respect to $b$ vanishes if
\begin{eqnarray*}
\lefteqn{\left[(u+M)c(c-b)^2 - (v-M)a(a+b)^2 + ac(a+c)(a+2b-c)\right] \nonumber } \\
&& \qquad \qquad \qquad \times \left[3(1+M^2) - 2M(a-c) - ac \right] = 0
\end{eqnarray*}
holds.
If the second factor vanishes, (\ref{pandq}) implies
\begin{eqnarray*}
\lefteqn{ p+q=\frac{3(1+M^2) - 2M(a-c) + b(c-a-b)}{(a+b)(c-b)} \nonumber }\\
&& \hspace{10em} = \frac{ac + bc - ab - b^2}{(a+b)(c-b)} =1, \qquad
\end{eqnarray*}
which is the boundary case ${\cal A}_5$ treated above.
So, ignoring the second factor and multiplying by $ac$ we arrive at
\begin{equation}\label{LA66a}
(u+M)ac^2(c-b)^2 - (v-M)a^2c(a+b)^2 + a^2c^2(a+c)(a+2b-c) = 0.
\end{equation}
The first equation from (\ref{LA61a}) can be rewritten as
\begin{equation}\label{LA67a}
(u+M)c^2(c-b)(3a^2 + 2a(b+c) + bc) + (v-M)a^2(a+b)^2 c - a^2c^2(a+c)^2 = 0.
\end{equation}
Adding up (\ref{LA66a}) and (\ref{LA67a}) and dividing the result by $c^2$ we obtain
\begin{equation}\label{LA68a}
(u+M)(c-b)(a+c)(3a+b) - 2a^2(a+c)(c-b) = 0
\end{equation}
and hence
\begin{eqnarray*}
b = \frac{2a^2}{u+M} - 3a.
\end{eqnarray*}
Substituting this into (\ref{LA66a}) with $c=a/\gamma$ and multiplying the result by $a^{-3}\gamma^3$ we get
\begin{eqnarray*}
\lefteqn{ 4 \left(\frac a{u+M} \right)^2 \left[\gamma(2\gamma +1)(u+M) -\gamma^3(v-M) \right] \nonumber } \\
&& + \frac a{u+M} \left[-(17\gamma^2 + 10\gamma + 1)(u+M) + 8\gamma^3 (v-M) \right] \\
&& + (3\gamma + 1)^2(u+M) -4\gamma^3(v-M) =0. \nonumber
\end{eqnarray*}
With the help of (\ref{LA63a}) this may be rewritten as
\begin{eqnarray}\label{LA611a}
\lefteqn{ 4 \left(\frac a{u+M} \right)^2 \left[(2\gamma^2 - 2\gamma - 1)(u+M) + 3\gamma^2(v-M) \right] \nonumber } \\
&& + \frac a{u+M} \left[(-17\gamma^2 + 14\gamma + 7)(u+M) - 24\gamma^2 (v-M) \right] \\
&& + 3(3\gamma + 1)(\gamma -1)(u+M) + 12\gamma^2(v-M) =0. \nonumber
\end{eqnarray}
One may verify that (\ref{LA611a}) can be factorized as follows
\begin{eqnarray*}
\lefteqn{ \left(\frac a{u+M} - 1 \right) \left\{\frac a{u+M}\left[(8\gamma^2 - 8\gamma - 4)(u+M) + 12\gamma^2(v-M) \right] \right. \nonumber } \\
&& \left. \qquad  - (9\gamma^2 - 6\gamma - 3)(u+M) - 12\gamma^2 (v-M) \right\} = 0.
\end{eqnarray*}
As $a=u+M$ is a boundary case, we conclude
\begin{eqnarray*}
a = \frac 34 (u+M) \frac{(\gamma - 1)(3\gamma + 1)(u+M) + 4 \gamma^2 (v-M)}{(2\gamma^2 - 2\gamma - 1)(u+M) + 3 \gamma^2 (v-M)}.
\end{eqnarray*}
Note that (\ref{LA63a}) may be reformulated as
\begin{eqnarray*}
(3\gamma +1)(u+M) = \gamma^2 (\gamma +3)(v-M)
\end{eqnarray*}
and that hence
\begin{eqnarray*}
a = \frac 34 (u+M) \frac{(\gamma + 1)^2 (v-M)}{(\gamma + 1)^2 (v-M) -(u+v)}.
\end{eqnarray*}
holds.
Consequently also $b, c,$ and the function to be maximized itself, as given in (\ref{functiontomaximizea}), can be expressed in terms of $M$ and $\gamma$.
As $\gamma$ is a complicated function of $M$, we shall write $M$ in terms of $\gamma$.
To this end we rewrite (\ref{LA63a}) as
\begin{eqnarray*}
\frac {u+M}{v-M} = \frac{\gamma^2 (\gamma +3)}{3\gamma + 1}
\end{eqnarray*}
and notice that this implies
\begin{equation}\label{LA617a}
M = \frac{\gamma^2(\gamma +3)v - (3\gamma +1)u}{(\gamma +1)^3}
\end{equation}
and hence
\begin{eqnarray}\label{LA618a}
&& a = \frac 38 \frac{\gamma(\gamma +3)(3\gamma +1)}{(\gamma +1)^3} (u+v), \nonumber \\
&& b = \frac 9{32} \frac{(1-\gamma)(\gamma +3)(3\gamma +1)}{(\gamma +1)^3} (u+v),  \\
&& c = \frac 38 \frac{(\gamma +3)(3\gamma +1)}{(\gamma +1)^3} (u+v). \nonumber
\end{eqnarray}

As $\gamma$ is positive, these values satisfy $-a<b<c$ and $-u \leq M \leq v$ as prescribed by $\cal A$.
Substituting them into (\ref{G0a})--(\ref{zetaeta}) we arrive at
\begin{eqnarray*}
\frac {256(\gamma +1)^5}{27(\gamma+3)^3(3\gamma +1)^3 (u+v)^2}\left[(\gamma +3)(\zeta + \eta c) + (3\gamma +1)(\zeta - \eta a) \right].
\end{eqnarray*}
Eliminating $\zeta$ and $\eta$ from this expression we obtain
\begin{eqnarray}\label{LA620a}
\lefteqn{ \left[9(\gamma+3)^3(3\gamma +1)^3 (u+v)^2\right]^{-1} \nonumber}\\
&& \left[ 1024 (\gamma +1)^6 + u^2(3\gamma + 1)^3 (9\gamma^3 + 33 \gamma^2 + 99 \gamma + 115) \right. \\
&& \ + v^2 (\gamma +3)^3(115 \gamma^3 + 99 \gamma^2 + 33\gamma + 9) \nonumber \\
&& \left. \ -2uv (\gamma+3)(3\gamma +1)(111 \gamma^4 + 276\gamma^3 + 250 \gamma^2 + 276 \gamma + 111) \right], \nonumber
\end{eqnarray}
which we shall denote as $\psi(\gamma; u,v)$.
At $\gamma =1$ the expression from (\ref{LA620a}) equals the first bound from (\ref{boundgeneralGauss}), i.e.,
\begin{eqnarray*}
\psi(1;u,v) = \frac 49\ \frac{4 + (u-v)^2}{(u+v)^2}.
\end{eqnarray*}
We shall prove
\begin{equation}\label{LA622a}
\psi(1;u,v) - \psi(\gamma;u,v) \geq 0,\quad \gamma >0, \quad \sqrt{3} \leq v \leq u.
\end{equation}
The function
\begin{eqnarray*}
\lefteqn{ \chi(\gamma; u,v) = 9(u+v)^2 (\gamma +3)^3 (3\gamma +1)^3 [\psi(1;u,v) - \psi(\gamma;u,v)] \nonumber } \\
&& = - 4(\gamma -1)^2 [148 \gamma^4 + 752 \gamma^3 + 1272 \gamma^2 + 752 \gamma + 148] \\
&& \quad + u^2 [ 189 \gamma^6 + 3186 \gamma^5 + 12051 \gamma^4 + 18556 \gamma^3 + 11667 \gamma^2 + 3186 \gamma + 317 ] \nonumber \\
&& \quad + v^2 [ 317 \gamma^6 + 3186 \gamma^5 + 11667 \gamma^4 + 18556 \gamma^3 + 12051 \gamma^2 + 3186 \gamma + 189 ] \nonumber \\
&& \quad - 2uv [ 99 \gamma^6 + 2382 \gamma^5 + 11853 \gamma^4 + 20484 \gamma^3 + 11853 \gamma^2 + 2382 \gamma + 99 ] \nonumber
\end{eqnarray*}
is quadratic in $u$ with a positive coefficient for $u^2$ and attains its minimum in $u$ at
\begin{equation}\label{LA624a}
v \frac{99 \gamma^6 + 2382 \gamma^5 + 11853 \gamma^4 + 20484 \gamma^3 + 11853 \gamma^2 + 2382 \gamma + 99}
{ 189 \gamma^6 + 3186 \gamma^5 + 12051 \gamma^4 + 18556 \gamma^3 + 11667 \gamma^2 + 3186 \gamma + 317},
\end{equation}
which equals at most $v$ as the denominator minus the numerator in (\ref{LA624a}) equals
\begin{eqnarray*}
2(\gamma - 1)^2 (45 \gamma^4 + 492 \gamma^3 + 1038 \gamma^2 + 620 \gamma + 109) \geq 0.
\end{eqnarray*}
We see that
\begin{eqnarray*}
\lefteqn{ \chi(\gamma;v,v) = -4 (\gamma -1)^2 [148 \gamma^4 + 752 \gamma^3 + 1272 \gamma^2 + 752 \gamma + 148] \nonumber } \\
&& + 4 v^2 (\gamma -1)^2 [77 \gamma^4 + 556 \gamma^3 + 1038 \gamma^2 + 556 \gamma + 77]
\end{eqnarray*}
is nonnegative for $\gamma >0$, if
\begin{eqnarray*}
v^2 \geq \sup_{\gamma >0} \frac {148 \gamma^4 + 752 \gamma^3 + 1272 \gamma^2 + 752 \gamma + 148} {77 \gamma^4 + 556 \gamma^3 + 1038 \gamma^2 + 556 \gamma + 77} = 2
\end{eqnarray*}
holds, which is the case in view of the assumption $v \geq \sqrt{3}$.
Consequently, $\chi(\gamma; u,v) \geq 0$ and (\ref{LA622a}) hold.
We have proved that the first bound from (\ref{boundgeneralGauss}) is valid for $\sqrt{3} \leq v \leq u$.

Substituting $\gamma=1$ into (\ref{LA617a}) and (\ref{LA618a}) we arrive at the random variables $Z$ and $Y$ defined in Theorem \ref{generalGaussa}.
However, $Y$ from this Theorem is not well defined if $2+v(v-u)$ is negative, i.e., if $u > v + 2/v$ holds.
Put differently, for $u > v + 2/v$ and $\gamma=1$ the point $(a,b,c,M)$ from (\ref{LA617a}) and (\ref{LA618a}) is not contained in $\cal A$.

Therefore, we have to consider the case $\sqrt{3} \leq v,\ v+2/v < u$ separately, as we do next.
Lengthy computations show that the derivative of $\psi(\gamma;u,v)$ from (\ref{LA620a}) with respect to $\gamma$ vanishes if and only if
\begin{eqnarray}\label{LA628a}
\lefteqn{ (\gamma^2 -1)\left[ 256(\gamma+1)^4 + (3\gamma +1)^4 u^2 + (\gamma+3)^4 v^2 \right. } \\
&& \hspace{10em} \left. -2(\gamma+3)(3\gamma+1)(29 \gamma^2 + 54 \gamma + 29)uv \right] =0 \nonumber
\end{eqnarray}
holds.
As $\gamma=1$ corresponds with a point outside $\cal A$ in the present situation, we restrict attention to those $\gamma$ for which the second factor vanishes.
So, we may assume
\begin{eqnarray}\label{LA628a}
\lefteqn{ - (3\gamma +1)^4 u^2 = 256(\gamma+1)^4 + (\gamma+3)^4 v^2 } \\
&& \hspace{10em} - 2(\gamma+3)(3\gamma+1)(29 \gamma^2 + 54 \gamma + 29)uv = 0. \nonumber
\end{eqnarray}
Using (\ref{pqpositive}), (\ref{LA617a}), and (\ref{LA618a}) we see that in terms of $\gamma$ the nonnegativity of $p$ is equivalent to
\begin{eqnarray*}
\lefteqn{ 256(\gamma+1)^6 -(\gamma -1)(3\gamma +1)^4 u^2 \nonumber } \\
&& + (\gamma+3)^2 \left( 112 \gamma^4 + 207 \gamma^3 + 139 \gamma^2 + 45 \gamma + 9 \right) v^2 \\
&& -2(\gamma+3)(3\gamma+1) \left( 24 \gamma^4 + 43 \gamma^3 + 79 \gamma^2 + 81 \gamma + 29 \right)uv \geq 0, \nonumber
\end{eqnarray*}
which by (\ref{LA628a}) becomes the inequality
\begin{eqnarray*}
\lefteqn{ 256 \gamma (\gamma+1)^4 (\gamma +3) +16 \gamma  (\gamma+1) (\gamma+3)^2 \left( 7 \gamma^2 + 6 \gamma + 3 \right) v^2 } \\
&& -16 \gamma (\gamma +1)(\gamma+3)(3\gamma+1) \left( 3 \gamma^2 + 6 \gamma + 7 \right)uv \geq 0. \nonumber
\end{eqnarray*}
In view of $u>v+2/v$ this implies
\begin{eqnarray*}
-32 \gamma (\gamma +1) (\gamma +3) (\gamma -1)^3 > 0
\end{eqnarray*}
and we conclude that for a $\gamma$ satisfying (\ref{LA628a}) $0 \leq \gamma <1$ holds.
Since $\gamma = \gamma_M$ from (\ref{gammaa}) is increasing in $(u+M)/(v-M)$ and equals 1 for $(u+M)/(v-M)=1$, we have $(u+M)/(v-M) < 1$ and hence
\begin{eqnarray*}
v + \frac 2v < u < v - 2M \leq v + 2\sqrt{3}.
\end{eqnarray*}
Applying these inequalities to (\ref{LA628a}) with the first factor removed we arrive at
\begin{eqnarray*}
\lefteqn{ 0 < 16(\gamma -1) \left( 55 \gamma^3 + 87 \gamma^2 + 45 \gamma + 5 \right) } \\
&& - 4 \left( 23 \gamma^4 + 196 \gamma^3 + 330 \gamma^2 + 196 \gamma + 23 \right) v^2 + 4 \sqrt{3} (3 \gamma +1)^4 v \nonumber
\end{eqnarray*}
and hence in view of $\sqrt{3} \leq v$ and $\gamma < 1$ at
\begin{eqnarray*}
\lefteqn{ 0 < 16(\gamma -1) \left( 55 \gamma^3 + 87 \gamma^2 + 45 \gamma + 5 \right) \nonumber } \\
&& \quad + 8 \left( 29 \gamma^4 -44 \gamma^3 - 138 \gamma^2 - 92 \gamma - 11 \right) v^2 \nonumber \\
&& < 16(\gamma -1) \left( 55 \gamma^3 + 87 \gamma^2 + 45 \gamma + 5 \right) \\
&& \quad - 8 \left( 15 \gamma^3 + 138 \gamma^2 + 92 \gamma + 11 \right) v^2 < 0. \nonumber
\end{eqnarray*}
This contradiction shows that the stationary points corresponding to the second factor in (\ref{LA628a}) do not belong to $\cal A$.
It follows that for $v + 2/v < u$ there are no stationary points within the interior of $\cal A$ and hence the supremum in (\ref{maximumoverA}) is attained at the boundary $\partial{\cal A}$ of $\cal A$.
Consequently, (\ref{maximumoverboundaryA}) completes the proof, as computation shows that the random variables $Z$ mentioned in the Theorem attain the bounds.
\hfill$\Box$ \\

{}


\begin{thebibliography}{}

\bibitem[{Bickel and Krieger (199)}]{Bickel}
\textsc{Bickel, P.J. and Krieger, A.M.} (1992).
\newblock Extensions of Chebychev's inequality with applications.
\newblock \textit{Probab. Math. Statist.} \textbf{13}, 293--310.

\bibitem[{Billingsley (1995)}]{Billingsley}
\textsc{Billingsley, P.} (1995).
\newblock \textit{Probability and Measure}.
\newblock Wiley, New York.

\bibitem[{Bienaym\'e (1853)}]{Bienayme}
\textsc{Bienaym\'e, I.J.} (1853).
\newblock Consid\'erations \`a l'appui de la decouverte de Laplace sur la loi de probabilit\'e dans la m\'ethode des moindres carr\'es.
\newblock \textit{C.R. hebd. S\'eance Acad. Sci. Paris} \textbf{37}, 309--24.
\newblock Reprinted (1867) \textit{Journal de math\'ematiques pures et appliqu\'es (2)} \textbf{12}, 158--176.

\bibitem[{Camp (1922)}]{Camp}
\textsc{Camp, B.H.} (1922).
\newblock A new generalization of Tchebycheff's statistical inequality.
\newblock \textit{Bulletin of the American Mathematical Society} \textbf{28}, 427--432.

\bibitem[{Cantelli (1928)}]{Cantelli}
\textsc{Cantelli, F.P.} (1928).
\newblock Sui confini della probabilit\`a.
\newblock \textit{Atti del Congresso Internazionale dei Matematici} \textbf{6}, 47--59.

\bibitem[{Chebyshev (1867)}]{Chebyshev}
\textsc{Chebyshev, P.} (1867).
\newblock (P. Tch\'ebychef) Des valeurs moyennes.
\newblock \textit{Journal de math\'ematiques pures et appliqu\'es (2)} \textbf{12}, 177--184.

\bibitem[{Clarkson et al. (2009)}]{Clarkson}
\textsc{Clarkson, E., Denny, J.L., and Shepp, L.} (2009).
\newblock ROC and the bounds on tail probabilities via theorems of Dubins and F. Riesz.
\newblock \textit{The Annals of Applied Probability} \textbf{19}, 467--476.

\bibitem[{Cram\'er (1946)}]{Cramer}
\textsc{Cram\'er, H.} (1946).
\newblock \textit{Mathematical Methods of Statistics}.
\newblock Princeton University Press, Princeton.

\bibitem[{DasGupta (2000)}]{DasGupta}
\textsc{DasGupta, A.} (2000).
\newblock Best constants in Chebyshev's inequalities with various applications.
\newblock \textit{Metrika} \textbf{51}, 185--200.

\bibitem[{Dharmadhikari and Joag-Dev (1988)}]{Dharmadhikari}
\textsc{Dharmadhikari, S.W. and Joag-dev, K.} (1988).
\newblock \textit{Unimodality, convexity and applications}.
\newblock Academic Press, Boston.

\bibitem[{Feller (1966)}]{Feller}
\textsc{Feller, W.} (1966).
\newblock \textit{An Introduction to Probability Theory and Its Applications}.
\newblock Wiley, New York.

\bibitem[{Ferentinos (1982)}]{Ferentinos}
\textsc{Ferentinos, K.} (1982)
\newblock On Tcebycheff’s type inequalities.
\newblock \textit{Trabajos Estad´ıst. Investigaci´on Oper.} \textbf{33}, 125-–132.

\bibitem[{Gau\ss \, (1823)}]{Gauss}
\textsc{Gau\ss \,, C.F.} (1823).
\newblock Theoria combinationis observationum erroribus minimis obnoxiae, pars prior.
\newblock \textit{Commentationes Societatis Regiae Scientiarum Gottingensis Recentiores} \textbf{5}, 1--58.

\bibitem[{Heyde and Seneta (1972)}]{Heyde}
\textsc{Heyde, C.C. and Seneta, E.} (1972).
\newblock Studies in the history of probability and statistics. XXXI.
\newblock The simple branching process, a turning point test and a fundamental inequality: a historical note on I.J. Bienaym\'e.
\newblock \textit{Biometrika} \textbf{59}, 680--683.

\bibitem[{Hooghiemstra and Van Mieghem (2015)}]{Hooghiemstra}
\textsc{Hooghiemstra, G. and Van Mieghem, P.F.A.} (2015).
\newblock An inequality of Gauss.
\newblock \textit{Nieuw Archief voor Wiskunde} \textbf{5/16}, 123--126.

\bibitem[{Ion (2001)}]{Ion}
\textsc{Ion, R.A.} (2001).
\newblock \textit{Nonparametric Statistical Process Control}.
\newblock PhD thesis University of Amsterdam, Amsterdam.

\bibitem[{Jensen (1906)}]{Jensen}
\textsc{Jensen, J.L.W.V.} (1906).
\newblock Sur les fonctions convexes et les in\'egalit\'es entre les valeurs moyennes.
\newblock \textit{Acta Mathematica} {\bf 30}, 175-–193.

\bibitem[{Khintchine (1938)}]{Khintchine}
\textsc{Khintchine, A.Ya.} (1938).
\newblock On Unimodal Distributions,
\newblock \textit{Izv. Nauchno Issled. Inst.Mat.Mech. Temsk. Gos. Univ.} {\bf 2}, 1--7.

\bibitem[{Klaassen et al. (2000)}]{Klaassen}
\textsc{Klaassen, C.A.J, Mokveld, P.J., and van Es, A.J.} (2000).
\newblock Squared skewness minus kurtosis bounded by 186/125 for unimodal distributions.
\newblock \textit{Statistics $\&$ Probability Letters} \textbf{50}, 131--135.

\bibitem[{Meidell (1922)}]{Meidell}
\textsc{Meidell, M.B.} (1922).
\newblock Sur un probl\`eme du calcul des probabilit\'es et les statistique math\'ematiques.
\newblock \textit{ Comptes Rendus} {\bf 175}, 806--808.

\bibitem[{Pukelsheim (1994)}]{Pukelsheim}
\textsc{Pukelsheim, F.} (1994).
\newblock The three sigma rule.
\newblock \textit{The American Statist.} {\bf 48}, 88--91

\bibitem[{Rougier et al. (2013)}]{Rougier}
\textsc{Rougier, J., Goldstein, M., and House, L.} (2013).
\newblock Second-Order Exchangeability Analysis for Multimodel Ensembles.
\newblock \textit{Journal of the American Statistical Association} {\bf 108}, 852--863.

\bibitem[{Selberg (1940)}]{Selberg}
\textsc{Selberg, H.L.} (1940)
\newblock Zwei Ungleichungen zur Erg\"anzung des Tchebycheffschen Lemmas.
\newblock \textit{Skand. Aktuarietidskr.} 1940 (1940), 121--125.

\bibitem[{Sellke (1996)}]{Sellke1}
\textsc{Sellke, T.} (1996).
\newblock Generalized Gau\ss \,-Chebyshev's Inequalities for Unimodal Distributions.
\newblock \textit {Metrika} \textbf{43}, 107--121.

\bibitem[{Sellke and Sellke (1997)}]{Sellke2}
\textsc{Sellke, T.M. and Sellke S.H.} (1997).
\newblock Chebyshev Inequalities for Unimodal Distributions.
\newblock \textit{The American Statist.} \textbf{51}, 34--40.

\bibitem[{Semenikhin (2019)}]{Semenikhin}
\textsc{Semenikhin, K.V.} (2019).
\newblock Two-Sided Probability Bound for a Symmetric Unimodal Random Variable.
\newblock \textit{Automation and Remote Control.} \textbf{80}, 474-–489.

\bibitem[{Shewhart (1931)}]{Shewhart}
\textsc{Shewhart, W.A.} (1931).
\newblock \textit{Economic Control of Quality of Manufactured Product}.
\newblock Van Nostrand, Princeton.

\bibitem[{Stewart (1995)}]{Stewart}
\textsc{Stewart, G.W.} (1995).
\newblock \textit{Theory of the Combination of Observations Least Subject to Error : part one, part two, supplement}.
\newblock \textit{Classics in Applied Mathematics} \textbf{11}.
\newblock Society for Industrial and Applied Mathematics, Philadelphia.
\newblock Translation from Latin of \cite{Gauss}.

\bibitem[{Theil (1949)}]{Theil}
\textsc{Theil, H.} (1949).
\newblock Over de ongelijkheid van Camp en Meidell.
\newblock \textit{Statistica Neerlandica} \textbf{3}, 201--208.

\bibitem[{Van den Heuvel and Ion (2003)}]{Heuvel}
\textsc{Van den Heuvel, E.R. and Ion, R.A.} (2003).
\newblock A Note on Capability Indices and the Proportion of Nonconforming Items.
\newblock \textit{Quality Engineering} \textbf{15}, 425-437.

\bibitem[{Vyso$\check{\rm c}$anski$\breve{\rm i}$ and Petunin (1980)}]{Vysochanskii1}
\textsc{Vyso$\check{\rm c}$anski$\breve{\rm i}$, D.F. and Petunin, Ju.I.} (1980).
\newblock Justification of the $3\sigma$-rule for unimodal distributions.
\newblock \textit{Theory of Probability and Mathematical Statistics} \textbf{21}, 25--36.

\bibitem[{Vysochanski$\breve{\rm i}$ and Petunin (1983)}]{Vysochanskii2}
\textsc{Vysochanski$\breve{\rm i}$, D.F. and Petunin, Yu.I.} (1983).
\newblock A remark on the paper "Justification of the $3\sigma$-rule for unimodal distributions".
\newblock \textit{Theory of Probability and Mathematical Statistics} \textbf{27}, 27--29.

\end{thebibliography}
\end{document}